\documentclass[twoside,twocolumn,fleqn,english,9pt,intoc,bibliography=totoc,index=totoc,BCOR10mm,captions=tableheading]{article}
\usepackage{lmodern}

\usepackage[T1]{fontenc}
\usepackage[utf8]{inputenc}
\usepackage[a4paper]{geometry}
\usepackage{babel}
\usepackage{float}
\usepackage{amsmath}
\usepackage{amssymb}
\usepackage{graphicx}
\usepackage{wasysym}
\usepackage{esint}
\usepackage{nomencl}

\providecommand{\makenomenclature}{\makeglossary}
\makenomenclature
\usepackage[unicode=true,
 bookmarks=true,bookmarksnumbered=true,bookmarksopen=true,bookmarksopenlevel=1,
 breaklinks=false,pdfborder={0 0 0},pdfborderstyle={},backref=false,colorlinks=false]
 {hyperref}
\hypersetup{pdftitle={On the Number of Conjugate Classes of Derangements},
 pdfauthor={Li, Wenwei},
 pdfsubject={Math.},
 pdfkeywords={Curve Fitting, Enumeration, Conjugate class, Derangements, Latin rectangles, Restricted Partition number, Estimation formula, Functional approximation, Accuracy},
 pdfpagelayout=OneColumn,pdfnewwindow=true,pdfstartview=XYZ,plainpages=false}

\makeatletter

\newcommand{\lyxdot}{.}

\newcommand{\lyxaddress}[1]{
	\par {\raggedright #1
	\vspace{1.4em}
	\noindent\par}
}

\usepackage[numbers,sort]{natbib} 

\usepackage{flushend,cuted}    

\usepackage{exscale}	
\usepackage{relsize}
\usepackage{mathcomp}  

\AtBeginDocument{
  
}

\makeatother

\begin{document}
\title{\textbf{On the Number of Conjugate Classes of Derangements}\textsf{}}

\maketitle
\begin{strip}
\begin{center} 
\author{\emph{Wen-Wei Li $^{1,2}$, $\quad\quad$ $\quad$ Zhong-Lin Cheng
$^{1}$ $\quad\quad$ }and\emph{ $\quad\quad$ Jia-Bao Liu $^{3}$}}

\lyxaddress{\begin{center}
\emph{liwenwei@ustc.edu}, \emph{$\quad\quad$ chzhlin}1234\emph{@}163.\emph{com},
\emph{$\quad\quad$} \emph{liujiabaoad}@163.\emph{com}$\ $ 
\par\end{center}}

\end{center}

\lyxaddress{1. School of Information and Mathematics, Anhui International Studies
University, Hefei, China, 231201\\
2. School of Mathematical Science, University of Science and Technology
of China, Hefei, China, 230026\\
3. School of Mathematics and Physics, Anhui Jianzhu University, Hefei,
China, 230601  }

\lyxaddress{Correspondence should be addressed to \emph{Wen-Wei Li}, \emph{liwenwei@ustc.edu} }

Published in \emph{Journal of Mathematics}, Volume 2021, Article ID
6023081, 20 pages, \\
https://doi.org/10.1155/2021/6023081
\begin{abstract}
The number of conjugate classes of derangements of order $n$ is the
same as the number $h(n)$ of the restricted partitions with every
portion greater than $1$. It is also equal to the number of isotopy
classes of $2\times n$ Latin rectangles. Sometimes the exact value
is necessary, while sometimes we need the approximation value. In
this paper, a recursion formula of $h(n)$ will be obtained, also
will some elementary approximation formulae with high accuracy for
$h(n)$ be presented. Although we may obtain the value of $h(n)$
in some computer algebra system, it is still meaningful to find an
efficient way to calculate the approximate value, especially in engineering,
since most people are familiar with neither programming nor CAS software.
This paper is mainly for the readers who need a simple and practical
formula to obtain the approximate value (without writing a program)
with more accuracy, such as to compute the value in an pocket science
calculator without programming function. Some methods used here can
also be applied to find the fitting functions for some types of 
data obtained in experiments.

\vspace{0.5cm}


\textbf{}

\textbf{AMS2020 Subject Classification:} 65D10, 05A17, 11P81, 

\vspace*{4cm}
\end{abstract}
\end{strip} 
\maketitle

\newpage{}

 \flushend 

\section{Introduction \label{sec:intro}}

Below $n$ is a positive integer greater than 1. 

On some occasions, it is necessary to know  the number of conjugate
classes of derangements. 

When generating the representatives of all the isotopy classes of
Latin rectangles of order $n$ by some method, we need to know the
number of the isotopy classes of $2\times n$ Latin rectangles for
verification. In some cases, we need the approximate value in  a
simple and efficient  method. (When writing a C program to generate
the representatives of all the isotopy classes of Latin rectangles
of order $n$, we need to prepare some space in the memory module
(RAM) to store the cycle structures of derangements so as to make
the program more efficient, otherwise, we have to allocate memory
dynamically, which will cost more time in memory addressing when writing
and reading data frequently in the particular position in the memory
module. So we need to know the number of the isotopy classes of $2\times n$
Latin rectangles for verification. )

 Let $\mathrm{S}_{n}$ be the symmetry group of the set $\mathrm{X}$
= \{1, 2, $\cdots$, $n$\}, i.e., the set (together with the operation
of combination) of the bijections from $\mathrm{X}$ to itself. An
element $\sigma$ in the symmetry group $\mathrm{S}_{n}$ is also
called a \emph{permutation} (of order $n$).   If $\sigma\in\mathrm{S}_{n}$,
$\sigma(i)\neq i$ ($\forall i\in\mathrm{X}$), $\sigma$ will be
called a \emph{derangement}\index{derangement} of order $n$.  If
a permutation $\sigma$ transforms any element in $\mathrm{X}$ to
a distinct element, then the sequence $\left[\sigma(1),\sigma(2),\cdots,\sigma(n)\right]$
will also be called a \emph{derangement}. The number of derangements
of order $n$ is denoted by $D_{n}$ (or $!n$ in some literature).
It is mentioned in nearly every combinatorics textbook that, 
\begin{align*}
D_{n} & =(n-1)\left(D_{n-1}+D_{n-2}\right)\\
 & =n!\sum_{i=0}^{n}\frac{(-1)^{i}}{i!}\doteq\left\lfloor \frac{n!}{\mathrm{e}}+\frac{1}{2}\right\rfloor ,\quad n\geqslant1.
\end{align*}
 Here $\left\lfloor x\right\rfloor $ is the floor function,  which
stands for the maximum integer that is less than or equal to the real
$x$.

For $x,y\in\mathrm{S}_{n}$, if $\exists z\in\mathrm{S}_{n}$, s.t.
$x=zyz^{-1}$, then $x$ and $y$ will be called \emph{conjugate,}
\index{conjugate} \index{permutationally similar} and $y$ is called
the \emph{conjugation }\index{conjugation} of $x$. Of course the
conjugacy relation is an equivalence relation. So the set of derangements
of order $n$ can be divided into some conjugate classes. This paper
is mainly concerned on the number of conjugate classes of derangements
of order $n$. The main method is  similar to that described in reference
\cite{liwenwei2016-Estmn-pn-arXiv}.

A matrix of size $k\times n$ ($1\leqslant k\leqslant n-1$) with
every row being a reordering of a fixed set of $n$ elements and every
column being a part of a reordering of the same set of $n$ elements,
is called a \emph{Latin rectangle}\index{Latin rectangle}. Usually,
the set of the $n$ elements is assumed to be \{$\,$1, 2, 3, $\cdots$,
$n$$\,$\}. (in some literatures, the members in a Latin rectangle
is assumed in the set \{$\,$0, 1, 2, $\cdots$, $n-1$$\,$\}.) 

 A $2\times n$ Latin rectangle with the first row in increasing
order could be considered as a derangement. An isotopy class of $2\times n$
Latin rectangles will correspond to a unique conjugate class of derangements.
 So  the number of isotopy classes of $2\times n$ Latin rectangles
is the same as the number of conjugate classes of derangements of
order $n$.

All the members in a conjugate class of derangements in $\mathrm{S}_{n}$
share the same cycle structure. Here we define the \emph{cycle structure}
of a derangement as the sequence in non-decreasing order of the lengths
(with duplicate entries) of all the cycles in the cycle decomposition
of the derangement. A cycle structure of a derangement of order $n$
could be considered as an integer solution of the equation 
\begin{equation}
s_{1}+s_{2}+\cdots+s_{q}=n,\quad\ (2\leqslant s_{1}\leqslant s_{2}\leqslant\cdots\leqslant s_{q}),\label{eq:Ptt-Rstrct-2}
\end{equation}
 where $s_{1}$, $s_{2}$, $\cdots$, $s_{q}$ are unknowns.

 For a fixed $q$, designate the number of integer solutions of the
equation \eqref{eq:Ptt-Rstrct-2} as $H_{q}(n)$, \label{Sym:Hq(n)}
\nomenclature[Hq(n)]{$H_q(n)$}{the number of integer solutions of the equation $s_{1} + s_{2} +\cdots +s_{q}$ =$n$  (2 $\leqslant$ $ s_{1}  \leqslant s_{2}$ $ \leqslant$  $\cdots$ $ \leqslant  s_{q}$) for a fixed positive integer $q$.  $H_q(n)>0$ when $q$  is less than$ \left\lfloor \frac{n}{2}\right\rfloor +1$. \pageref{Sym:Hq(n)}}
where $q$ is less than $\left\lfloor \frac{n}{2}\right\rfloor +1$
(otherwise $H_{q}(n)$ is defined by 0), and denote $h(n)$ the number
of all the integer solutions of Equation \eqref{eq:Ptt-Rstrct-2}
for all possible $q$, \label{Sym:h(n)} \nomenclature[h(n)]{$h(n)$}{the number of integer solutions of the equation $s_{1} + s_{2} +\cdots +s_{q}$ =$n$  (2 $\leqslant$ $ s_{1}  \leqslant s_{2}$ $ \leqslant$  $\cdots$ $ \leqslant  s_{q}$).  \pageref{Sym:h(n)}}
i.e., 
\[
h(n)=\underset{{\scriptstyle {\scriptstyle q=1}}}{\overset{{\scriptstyle \left\lfloor \frac{n}{2}\right\rfloor }}{\sum}}H_{q}(n).
\]
So the number of conjugate classes of derangements of order $n$ is
$h(n)$. Since $h(n)$ is the number of a type of restricted partitions,
it is tightly connected with the partition number $p(n)$. 

Following the notation of \cite{Marshall1958Survey}, denote by $P_{q}(n)$
\label{Sym:Pq(n)} \nomenclature[Pq(n)]{$P_q(n)$}{the number of integer solutions of the equation $s_{1} + s_{2} +\cdots +s_{q}$ =$n$  ( 1 $\leqslant$ $ s_{1}  \leqslant s_{2}$ $ \leqslant$  $\cdots$ $ \leqslant  s_{q}$) for a fixed positive integer $q$.  $P_q(n)>0$ when $q$  is less than$ \left\lfloor \frac{n}{2}\right\rfloor +1$. \pageref{Sym:Pq(n)}}
the number of integer solutions of equation 
\begin{equation}
s_{1}+s_{2}+\cdots+s_{q}=n,\quad\ (1\leqslant s_{1}\leqslant s_{2}\leqslant\cdots\leqslant s_{q})\label{eq:Ptt-Genrl}
\end{equation}
 for a fixed $q$, where $1\leqslant q\leqslant n$, and by $p(n)$
\label{Sym:p(n)} \nomenclature[p(n)]{$p(n)$}{the number of integer solutions of the equation $s_{1} + s_{2} +\cdots +s_{q}$ =$n$  (1 $\leqslant$ $ s_{1}  \leqslant s_{2}$ $ \leqslant$  $\cdots$ $ \leqslant  s_{q}$).  \pageref{Sym:p(n)}}
the number of all the (unrestricted) partitions of $n$. It is clear
that \footnote{$\ $ In a lot of articles, $p(n,q)$ is used in stead of $P_{q}(n)$,
but in some other literatures, $p(n,q)$ stands for some other number.} 
\begin{equation}
p(n)=\underset{_{q=1}}{\overset{_{n}}{\sum}}P_{q}(n).
\end{equation}

There is a brief introduction of the important results on the partition
number (or partition function) $p(n)$ and $P_{q}(n)$ in reference
\cite{Marshall1958Survey}, such as the recursion formula of $p(n)$
and $P_{q}(n)$. More information about the partition number $p(n)$
may be found in reference \cite{Eric1999PtFuncP}. There is a list
of some important papers and book chapters on the partition number
in \cite{SloaneOEIS-A000041-pn} (including the ``LINKS'' and ``REFERENCES
'') and \cite{Anonymous2007Bib-Partition}.  Reference \cite{liwenwei2016-Estmn-pn-arXiv}
presented some estimation formulae with high accuracy for $p(n)$,
which are revised from the Hardy-Ramanujan’s asymptotic formula. 

There are also a lot of literatures on the number of some types of
restricted partitions of $n$  (such as \cite{Newman1959W-RestrPart},
\cite{Lahiri1969SRPFCM3}, \cite{Lahiri1969SRPFCM5}, \cite{Lahiri1969SRPFCM7},
 \cite{Rafael2009ResPtEleMeth}, \cite{Brandt2012ARRLCPRPFABV},
\cite{Amanda2012-l-adic-Properties-Elsevier}, \cite{Eva2014-l-adic-Properties-Springer})
or on the congruence properties of (restricted) partition function
(such as \cite{Thanigasalam1974CPCRP}, \cite{Dean1975Note-CPCRP},
 \cite{Matt2002congruencesof}, \cite{Brandt2005ConPro-p}, \cite{Brandt2007OCPCVP},
\cite{Bruce2007Ramanujan'sCongruences}, \cite{George2012RUMPTF},
\cite{Brandt2013GenConPResPtFun}).  

In \cite{OEISwiki2011RestrPt}, we can find many cases of Restricted
Partitions (some of them are introduced in \cite{NIST2015IPRNPS},
\cite{NIST2015IPOR} or \cite{Frank2010NHMF}).  One class are concerned
on the restriction of the sizes of portions, such as portions restricted
to Fibonacci numbers, powers (of 2 or 3), unit, primes, non-primes,
composites or non-composites; another class are related to the restriction
of the number of portions, such as the cases that the number of parts
will not exceed $k$; the third class are about the restrictions for
both, for example, the cases that the number of parts is restricted
while the parts restricted to powers  or primes. But the author has
not found too much information on the number $h(n)$, especially on
the approximate calculation, although we can find a lot of information
on other restricted partition numbers.

The basic method of function fitting (curve fitting) could be found
in any textbook of numerical analysis, such as \cite{Anne2012NMDACIA}
and \cite{David2002NMMSC}. Some tricks used here may be found in
some books of data analysis such as \cite{William1998DataAnaly}.
They may be helpful for understanding the methods described in Section
\ref{sec:Estimate-h(n)}. 

Section \ref{sec:Some-Formulae-for} will deduce the recursion formula
for $h(n)$ and will show the relation of $h(n)$ and $p(n)$. Subsection
\ref{subsec:Asymptotic-Formula} will deduce the asymptotic formula
of $h(n)$ from Hardy and Ramanujan's Asymptotic formula of $p(n)$
(mentioned in \cite{liwenwei2016-Estmn-pn-arXiv}). This new asymptotic
formula $I_{\mathrm{g}}(n)$ coincides with Ingham's result (refer
\cite{Ingham1941TaubThm4Pt} and \cite{Daniel2006EleDerAsyPtFunc}).
By bringing in two parameters $C_{1}(n)$ and $C_{2}(n)$ in the new
asymptotic formula $I_{\mathrm{g}}(n)$, we have reached several estimation
formulae for $h(n)$ with high accuracy in subsection \ref{subsec:Method-B},
using the same  ideas described in \cite{liwenwei2016-Estmn-pn-arXiv}.
By fitting $h(n)-I_{\mathrm{g}}(n)$, we have another two estimation
formulae for $h(n)$ in subsection \ref{subsec:Method-C}. When $n<100$,
we have a more accurate estimation formula for $h(n)$ in subsection
\ref{subsec:Estimate-less-100}. The relative errors of these estimation
formulae will be presented to shown  the accuracy.

\section{Some Formulae for $h(n)$ \label{sec:Some-Formulae-for}}

In this section, we will derive a recursive formula from the method
mentioned in reference \cite{Marshall1958Survey} (page 53\textasciitilde 55).

By definition, $h(k)$ = 0 when $k<2$, but here we assume that $h(k)$
= 0 when $k<0$, $h(0)=1$ and $h(1)=0$, for convenience.

It is mentioned in \cite{Marshall1958Survey} (page 52) that in 1942
Auluck gave an estimation of $P_{q}(n)$ by $P_{q}(n)$ $\approx\dfrac{1}{q!}\left(\begin{array}{c}
n-1\\
q-1
\end{array}\right)$ when $n$ is large enough. 

By the same method shown in reference \cite{Marshall1958Survey} (page
53, 57), we can obtain the generation function of $h(n)$: 
\begin{align}
G(x) & =\sum\limits _{n=0}^{\infty}h(n)x^{n}=\prod\limits _{i=2}^{\infty}\left(1-x^{i}\right)^{-1}\nonumber \\
 & =\frac{1}{\left(1-x^{2}\right)}\frac{1}{\left(1-x^{3}\right)}\frac{1}{\left(1-x^{4}\right)}\cdots\frac{1}{\left(1-x^{i}\right)}\cdots\cdots\label{eq:hn-gen-function}
\end{align}
and a formula 
\begin{equation}
h(n)=\frac{1}{2\pi i}\mathlarger{\oint}_{C}\frac{G(x)}{x^{n+1}}\mbox{d}x,
\end{equation}
 where $h(0)=1$, $h(1)=0$, and $C$ is a contour around the original
point. 

It is difficult to get a simple formula to count the solutions of
Equation \eqref{eq:Ptt-Rstrct-2} in general. But for a fixed integer
$q$, the number $H_{q}(n)$ of solutions is $0$ (when $q>\left\lfloor \frac{n}{2}\right\rfloor $)
or \\
 $\underset{s_{1}=2}{\overset{\left\lfloor \frac{n}{q}\right\rfloor }{\sum}}\underset{s_{2}=s_{1}}{\overset{\left\lfloor \frac{n-s_{1}}{q-1}\right\rfloor }{\sum}}\cdots\underset{s_{q-1}=s_{q-2}}{\overset{{\scriptscriptstyle \left\lfloor \frac{n-s_{1}-s_{2}\cdots-s_{q-2}}{2}\right\rfloor }}{\sum}}1$
\\
 = $\underset{s_{1}=2}{\overset{\left\lfloor \frac{n}{q}\right\rfloor }{\sum}}\underset{s_{2}=s_{1}}{\overset{\left\lfloor \frac{n-s_{1}}{q-1}\right\rfloor }{\sum}}\cdots\underset{{\scriptstyle s_{q-2}=s_{q-3}}}{\overset{\left\lfloor \footnotesize\frac{n-s_{1}-s_{2}\cdots-s_{q-3}}{3}\right\rfloor }{\sum}}$
\\
$\negthickspace$ \phantom{= =} $\left(\frac{n-s_{1}-s_{2}\cdots-s_{q-2}}{2}-s_{q-2}+1\right)$
\\
= $P_{q}(n-q)$ (when $q\leqslant\left\lfloor \frac{n}{2}\right\rfloor $).

 It is not difficult to find out that
\begin{equation}
h(n)=\underset{{\scriptstyle q=1}}{\overset{{\scriptscriptstyle \left\lfloor \frac{n}{2}\right\rfloor }}{\sum}}H_{q}(n)=\underset{{\scriptstyle q=1}}{\overset{{\scriptscriptstyle \left\lfloor \frac{n}{2}\right\rfloor }}{\sum}}P_{q}(n-q).\label{eq:hn-calculate}
\end{equation}
And there is a recursion for $P_{q}(n)$ in reference \cite{Marshall1958Survey}
(page 51) 
\begin{equation}
P_{q}(n)=\underset{{\scriptstyle j=1}}{\overset{{\scriptstyle t}}{\sum}}P_{j}(n-q),
\end{equation}
where $t=\min\{q,\,n-q\}$, so there is no difficulty to obtain the
values of $P_{q}(n)$ and $h(n)$ when $n$ is small.

For the value of $p(n)$ there is a recursion, 
\begin{align}
p(n) & =p(n-1)+p(n-2)-p(n-5)\nonumber \\
 & \ \ \ \ -p(n-7)+\cdots+\nonumber \\
 & \ \ \ \ (-1)^{k-1}p\left(n-\frac{3k^{2}\pm k}{2}\right)+\cdots\cdots\nonumber \\
 & =\sum\limits _{k=1}^{k_{1}}(-1)^{k-1}p\left(n-\frac{3k^{2}+k}{2}\right)+\nonumber \\
 & \ \ \ \ \sum\limits _{k=1}^{k_{2}}(-1)^{k-1}p\left(n-\frac{3k^{2}-k}{2}\right),\label{eq:pn-recursion}
\end{align}
where 
\begin{equation}
k_{1}=\left\lfloor \frac{\sqrt{24n+1}-1}{6}\right\rfloor ,\ k_{2}=\left\lfloor \frac{\sqrt{24n+1}+1}{6}\right\rfloor ,\label{eq:pn-recursion-k1-k2}
\end{equation}
 and assume that $p(0)=1$. (Refer \cite{Marshall1958Survey}, page
55). Here we assume that $p(x)$ = 0 when $x<0$.

We can obtain the same recursion for $h(n)$, 
\begin{align}
h(n) & =h(n-1)+h(n-2)-h(n-5)\nonumber \\
 & \ \ \ \ -h(n-7)+\cdots+\nonumber \\
 & \ \ \ \ (-1)^{k-1}h\left(n-\frac{3k^{2}\pm k}{2}\right)+\cdots\cdots\nonumber \\
 & =\sum\limits _{k=1}^{k_{1}}(-1)^{k-1}h\left(n-\frac{3k^{2}+k}{2}\right)+\nonumber \\
 & \ \ \ \ \sum\limits _{k=1}^{k_{2}}(-1)^{k-1}h\left(n-\frac{3k^{2}-k}{2}\right),\label{eq:hn-recursion}
\end{align}
where $k_{1}$ and $k_{2}$ are determined by Equation \eqref{eq:pn-recursion-k1-k2}
and assume that $h(0)=1$, $h(k)$ = 0 when $k<0$.

The proof of Equation \eqref{eq:hn-recursion} is easy to understand. 

By Equation \eqref{eq:hn-gen-function}, we have 
\begin{equation}
\left(\sum\limits _{n=0}^{\infty}h(n)x^{n}\right)\left(\prod\limits _{i=2}^{\infty}\left(1-x^{i}\right)\right)=1.\label{eq:hn-prop-1}
\end{equation}

 Since $F(x)=\sum\limits _{n=0}^{\infty}p(n)x^{n}=\prod\limits _{i=1}^{\infty}\left(1-x^{i}\right)^{-1}$,
where $p(0)=1$. \\
So $\left(\sum\limits _{n=0}^{\infty}p(n)x^{n}\right)\left(\prod\limits _{i=1}^{\infty}\left(1-x^{i}\right)\right)=1$,
or 
\begin{equation}
\left(\sum\limits _{n=0}^{\infty}p(n)x^{n}\right)(1-x)\left(\prod\limits _{i=2}^{\infty}\left(1-x^{i}\right)\right)=1.\label{eq:pn-prop-1}
\end{equation}

Compare Equation \eqref{eq:hn-prop-1} and Equation \eqref{eq:pn-prop-1},
we have 
\begin{align*}
\sum\limits _{n=0}^{\infty}h(n)x^{n} & =\left(\sum\limits _{n=0}^{\infty}p(n)x^{n}\right)(1-x)\\
 & =\sum\limits _{n=0}^{\infty}\bigl(p(n)-p(n-1)\bigl)x^{n},
\end{align*}

by assumption $h(k)=p(k)=0$ when $k<0$. Hence,
\begin{equation}
h(n)=p(n)-p(n-1),\ \ (n=0,1,2,\cdots).\label{eq:Rel-hn-pn}
\end{equation}

By Equation \eqref{eq:Rel-hn-pn}, 
\begin{align}
 & h(n)\nonumber \\
= & h(n-1)+h(n-2)-h(n-5)-h(n-7)+\nonumber \\
 & \ \ \cdots+(-1)^{k-1}h\left(n-\frac{3k^{2}\pm k}{2}\right)+\cdots\cdots\nonumber \\
= & \sum\limits _{k=1}^{k_{1}}(-1)^{k-1}h\left(n-\frac{3k^{2}+k}{2}\right)+\nonumber \\
 & \ \sum\limits _{k=1}^{k_{2}}(-1)^{k-1}h\left(n-\frac{3k^{2}-k}{2}\right).\label{eq:hn-recursion-V2}
\end{align}
 We can easily obtain the solutions of Equation \eqref{eq:Ptt-Rstrct-2}
by hand when $n$ < 13. By Equation \eqref{eq:hn-recursion}, we can
obtain the number $h(n)$ of solutions of Equation \eqref{eq:Ptt-Rstrct-2}
  by writing a small program in some Computer Algebra System (CAS)
softwares such as ``maple'', ``maxima'', ``axiom'' or some other
softwares likewise (be aware of that 0 is not a valid index value
in some software such like maple).

\begin{table}[!b]
\begin{centering}
\center\includegraphics[viewport=97bp 461bp 405bp 703bp,clip,scale=0.63]{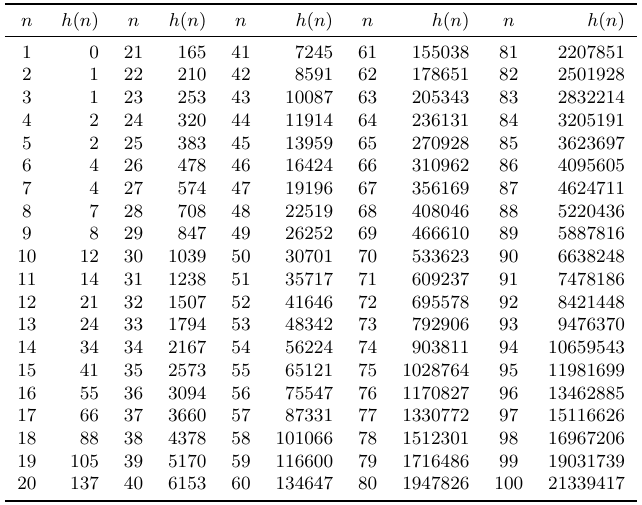}
\par\end{centering}
\caption{The value of $h(n)$ when $1\leqslant n\leqslant100$}
 \label{Table:N_hn}

\vspace{0.8cm}

\begin{centering}
\includegraphics[viewport=97bp 408bp 446bp 596bp,clip,scale=0.63]{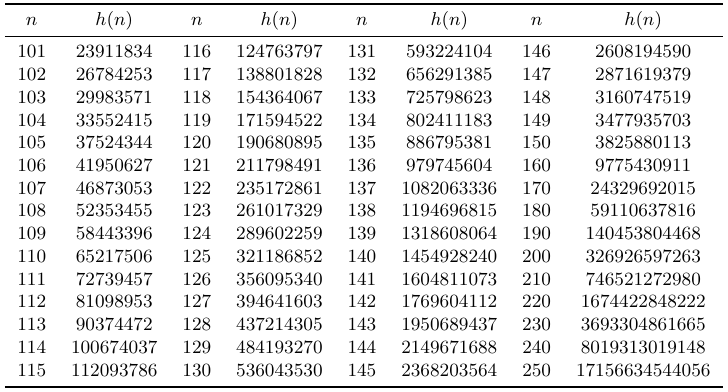}
\par\end{centering}
\begin{raggedright}
\caption{The value of $h(n)$ when $101\leqslant n\leqslant250$}
 \label{Table:N_hn-2}\vspace{0.8cm}
\par\end{raggedright}
\begin{centering}
\includegraphics[viewport=97bp 274bp 407bp 429bp,clip,scale=0.63]{tables-v2/Table-2_4_2_Pre_1_2_value_of_h-n_h-k-n_1_100}
\par\end{centering}
\caption{The number of solutions of Equation \eqref{eq:Ptt-Rstrct-2} for different\emph{
$q$}}
 \label{Table:N_Hqn}
\end{table}

The value of $h(n)$ when $n$ < 250 are shown on Table \ref{Table:N_hn}
(on page \pageref{Table:N_hn}) and Table \ref{Table:N_hn-2} (on
page \pageref{Table:N_hn-2}). Some value of $H_{q}(n)$ are shown
on Table \ref{Table:N_Hqn} (on page \pageref{Table:N_Hqn}).

Obviously, $h(n)$ < $p(n)$ holds by definition (when $n>1$). As
$p(n)$ grows much more slowly than exponential functions, i.e., for
any $r>1$, $p(n)<r^{n}$ will hold when $n$ is large enough, which
means we can not estimate $p(n)$ and $h(n)$ by an exponential function.
As $p(n)$ grows faster than any power of $n$, which means we can
not estimate $p(n)$ by a polynomial function. (refer \cite{Marshall1958Survey},
page 53) So, $h(n)$ can not be estimated by a polynomial function,
either. 

\section{The Estimation of $h(n)$ \label{sec:Estimate-h(n)}}

The recursion formula Equation \eqref{eq:hn-recursion} for $h(n)$
is not convenient in practical for a lot of people who do not want
to write programs. Sometimes we need the approximation value, such
as the cases mentioned in \cite{liwenwei2016-Estmn-pn-arXiv}, so
an estimation formula is necessary.

The figure of the data $\bigl(n,\,\ln(h(n))\bigr)$ ( $n$ = $60+20k$,
$k$ = 1, 2, $\cdots$, 397) are shown on  Figure  1 on page \pageref{Fig:2_4_24_(n-ln(h(n)))}.
  Here the data points are displayed by small hollow circles, and
the circles are very crowded that we may believe that the circles
themselves be a very thick curve if we notice only the right-hand
part. In this figure, the data points in the lower left part are sparse
(compared with the points in the upper right part), and we may find
some hollow circles easily. If there is a curve passes through these
hollow circles, we will notice it (as shown on Figure \ref{Fig:2_4_26_(n-ln(h(n)))-Ig1(n)}
on page \pageref{Fig:2_4_26_(n-ln(h(n)))-Ig1(n)}). But later in Figure
\ref{Fig:2_4_25_(n-C4(n))}, the circles distribute uniformly on a
curve, it will be difficult to distinguish the circles and the curve
passes through the centers of the them.

The author has not found a practical estimation formula with high
accuracy of the number $h(n)$ before. 

\begin{figure}
 \includegraphics[scale=0.25]{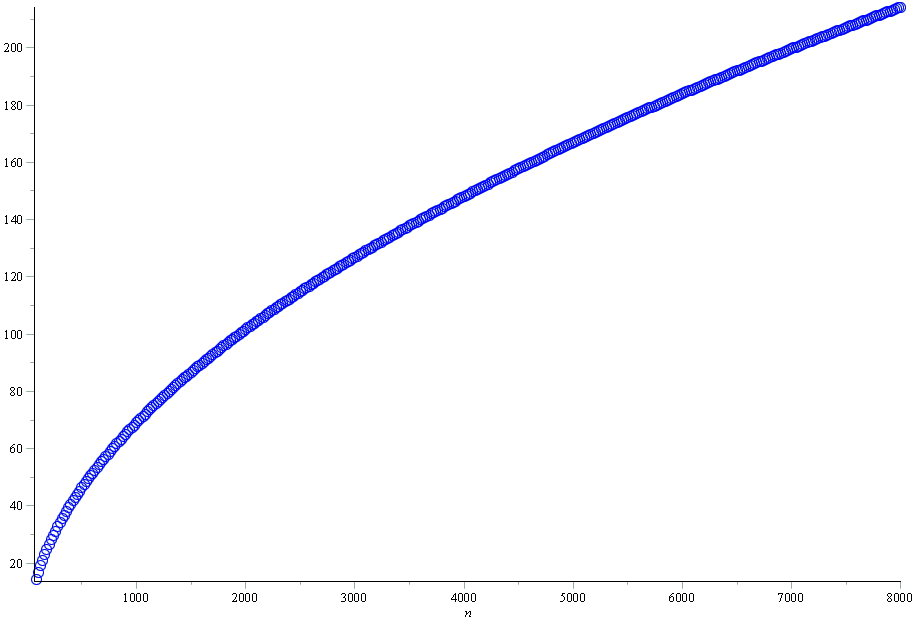}\label{Fig:2_4_24_(n-ln(h(n)))}

\caption{The graph of the data $\left(n,\ln h(n)\right)$}
\end{figure}

Actually, it is very difficult to find directly a simple function
to fit the data on Figure 1  with high accuracy. The main reason
is that the fitting functions obtained by the methods used frequently
could not reach satisfying accuracy.

Since we have several accurate estimation formula of $p(n)$ (refer
\cite{liwenwei2016-Estmn-pn-arXiv}), such as 
\[
R'_{\mathrm{h2}}(n)=\left\lfloor \frac{\exp\left(\sqrt{\frac{2}{3}}\pi\sqrt{n}\right)}{4\sqrt{3}\left(n+a_{2}\sqrt{n+c_{2}}+b_{2}\right)}+\frac{1}{2}\right\rfloor ,\ (n\geqslant80)
\]
 and 
\[
R'_{\mathrm{h0}}(n)=\left\lfloor \frac{\exp\left(\sqrt{\frac{2}{3}}\pi\sqrt{n}\right)}{4\sqrt{3}\left(n+C'_{2}(n)\right)}+\frac{1}{2}\right\rfloor ,\quad1\leqslant n\leqslant100,
\]
 where $a_{2}=0.4432884566$, $b_{2}=0.1325096085$, $c_{2}=0.274078$
and 
\begin{equation}
C'_{2}(n)=\begin{cases}
0.4527092482\times\sqrt{n+4.35278}-\\
\quad0.05498719946,\quad\ n=3,5,7,\cdots,99;\\
0.4412187317\times\sqrt{n-2.01699}+\\
\quad0.2102618735,\quad\quad n=4,6,8\cdots,100.
\end{cases}\label{eq:C2'(n)}
\end{equation}

By Equation \eqref{eq:Rel-hn-pn}, we can obtain $h(n)$ by 
\begin{equation}
h_{1}(n)=\begin{cases}
R'_{\mathrm{h0}}(n)-R'_{\mathrm{h0}}(n-1), & 2\leqslant n\leqslant80;\\
R'_{\mathrm{h2}}(n)-R'_{\mathrm{h2}}(n-1), & n>80.
\end{cases}\label{eq:h(n)-estimation-0}
\end{equation}
 and the error of this formula will not exceed twice of the error
of $R'_{\mathrm{h2}}(n)$ or $R'_{\mathrm{h0}}(n)$. Of course, this
formula is not as simple  as we want, but the accuracy is very good.

\subsection{Asymptotic Formula\label{subsec:Asymptotic-Formula}}

As $h(n)=p(n)-p(n-1)$, by \emph{Hardy-Ramanujan's asymptotic formula\index{Hardy-Ramanujan asymptotic formula}}
\[
p(n)\sim\frac{1}{4n\sqrt{3}}\exp\left(\sqrt{\frac{2}{3}}\pi\sqrt{n}\right)
\]
 (refer \cite{Ramanujan1918AsymFmlComAnal}, \cite{Pal1942ElemProofRHFml},
\cite{Donald1951EvalfConstHRFml}, \cite{Newman1962SimProofPtFml},
\cite{Eric1999PtFuncP}, \cite{NIST2015FNTANTUP}, \cite{liwenwei2016-Estmn-pn-arXiv}),
we assume that, when $n\gg1$, \\
$h(n)$ $\sim$ $\frac{\exp\left(\sqrt{\frac{2}{3}}\pi\sqrt{n}\right)}{4\sqrt{3}n}-\frac{\exp\left(\sqrt{\frac{2}{3}}\pi\sqrt{n-1}\right)}{4\sqrt{3}(n-1)}$.
So,

 $h(n)$ $\sim$ $\frac{\exp\left(\sqrt{\frac{2}{3}}\pi\sqrt{n-1}\right)}{4\sqrt{3}}\left(\frac{\exp\left(\pi\sqrt{\frac{2}{3}}\left(\sqrt{n}-\sqrt{n-1}\right)\right)}{n}-\frac{1}{n-1}\right)$

= $\frac{\exp\left(\sqrt{\frac{2}{3}}\pi\sqrt{n-1}\right)}{4\sqrt{3}}\left(\frac{\exp\left(\frac{\pi\sqrt{2/3}}{\sqrt{n}+\sqrt{n-1}}\right)}{n}-\frac{1}{n-1}\right)$

$\sim$ $\frac{\exp\left(\sqrt{\frac{2}{3}}\pi\sqrt{n}\right)}{4\sqrt{3}}\left(\frac{\exp\left(\frac{\pi\sqrt{2/3}}{2\sqrt{n}}\right)}{n}-\frac{1}{n-1}\right)$

$\sim$ $\frac{\exp\left(\sqrt{\frac{2}{3}}\pi\sqrt{n}\right)}{4\sqrt{3}}\left(\frac{1+\frac{\pi}{\sqrt{6n}}}{n}-\frac{1}{n-1}\right)$

$\sim$ $\frac{\exp\left(\sqrt{\frac{2}{3}}\pi\sqrt{n}\right)}{4\sqrt{3}}\left(\frac{\pi}{\sqrt{6n^{3}}}\right)$$=$
$\frac{\pi\exp\left(\sqrt{\frac{2}{3}}\pi\sqrt{n}\right)}{12\sqrt{2n^{3}}}.$

So, 
\begin{equation}
h(n)\sim\frac{\pi}{12\sqrt{2n^{3}}}\exp\left(\sqrt{\frac{2}{3}}\pi\sqrt{n}\right).\label{eq:h(n)-asymptotic}
\end{equation}

\begin{table}[h]

\begin{centering}
\includegraphics[viewport=98bp 479bp 437bp 666bp,clip,scale=0.67]{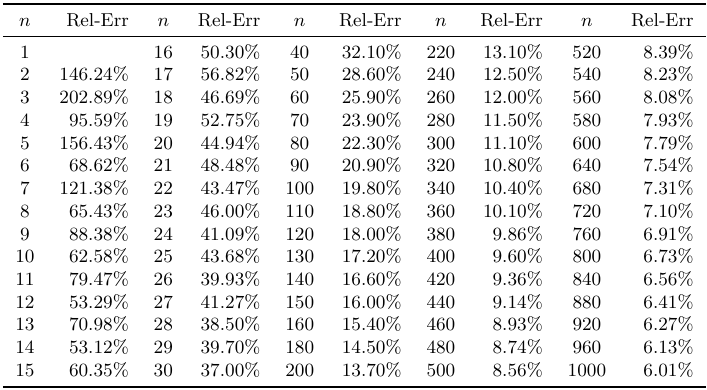}
\par\end{centering}
\caption{The relative error of $I_{\mathrm{g}}(n)$ to $h(n)$ when $n\leqslant1000$.}
 \label{Table:Rel-Err-h(n)-Ig(x)} \vspace{0.5cm}

\begin{centering}
\includegraphics[viewport=97bp 358bp 412bp 447bp,clip,scale=0.68]{tables-v2/Table-2_4_2_a-Rel-Err-Ig_n_-Ig_n_round}
\par\end{centering}
\caption{The relative error of $\left\lfloor I_{\mathrm{g}}(n)+\frac{1}{2}\right\rfloor $
to $h(n)$ when $n\leqslant30$.}
 \label{Table:Rel-Err-h(n)-Ig(n)-round}

\end{table}

In coincidence, half a year after the main results were obtained in
this paper, the author found an asymptotic formula
\begin{align}
P_{a,b}(n)\text{\ensuremath{\sim}} & \Gamma\left(\frac{b}{a}\right)\pi^{b/a-1}2^{-(3/2)-(b/2a)}3^{-(b/2a)}\nonumber \\
 & \quad a^{-(1/2)+(b/2a)}n^{-\frac{a+b}{2a}}\mathrm{e}^{\pi\sqrt{\frac{2n}{3a}}},\label{eq:P_a,b}
\end{align}

in \cite{Daniel2006EleDerAsyPtFunc}. When $a$ = 1, $b$ = 2, we
will have  
\begin{equation}
P_{1,2}(n)\text{\ensuremath{\sim}}\frac{\pi}{12\sqrt{2n^{3}}}\exp\left(\pi\sqrt{\frac{2}{3}n}\right),
\end{equation}
 which coincides with the asymptotic formula obtained here. 

The formula \eqref{eq:h(n)-asymptotic} will  be called the\emph{
Ingham-Meinardus asymptotic formula} \index{Ingham-Meinardus asymptotic formula}
in this paper, since Daniel mentioned in \cite{Daniel2006EleDerAsyPtFunc}
that this general asymptotic formula \eqref{eq:P_a,b} was first given
by A. E. Ingham in \cite{Ingham1941TaubThm4Pt} and the proof was
refined by G. Meinardus later in another two papers written in German.

Later in this paper, $\frac{\pi}{12\sqrt{2n^{3}}}\exp\left(\sqrt{\frac{2}{3}}\pi\sqrt{n}\right)$
will be denoted by $I_{\mathrm{g}}(n)$ for short. \label{Sym:Ig(n)}
\nomenclature[Ig(n)]{$I_{\mathrm{g}}(n)$}{The Ingham-Meinardus asymptotic formula. \pageref{Sym:Ig(n)}}

It is not satisfying to estimate $h(n)$ by $I_{\mathrm{g}}(n)$ when
$n$ is small. The relative error of $I_{\mathrm{g}}(n)$ to $h(n)$
is greater than 6\% as shown on Table \ref{Table:Rel-Err-h(n)-Ig(x)}
(on page \pageref{Table:Rel-Err-h(n)-Ig(x)}). The round approximation
\[
I'_{\mathrm{g}}(n)=\left\lfloor I_{\mathrm{g}}(n)+\frac{1}{2}\right\rfloor 
\]
 will not change the accuracy distinctly, as shown on Table \ref{Table:Rel-Err-h(n)-Ig(n)-round}
(on page \pageref{Table:Rel-Err-h(n)-Ig(n)-round}). So it is necessary
to modify the asymptotic formula  for better accuracy. 

\subsection{Method A: Modifying the exponent  \label{subsec:Modifying-the-exponent}}

In this subsection we consider fitting $h(n)$ by $I_{\mathrm{ga}}=\frac{\pi}{12\sqrt{2n^{3}}}\exp\left(\sqrt{\frac{2}{3}}\pi\sqrt{n+C_{1}(n)}\right)$,
or fitting \\
 $\left(n,\,\frac{3}{2\pi^{2}}\left(\ln\left(\frac{12\sqrt{2n^{3}}h(n)}{\pi}\right)\right)^{2}-n\right)$
( $n$ = $60+20k$, $k$ = 1, 2, $\cdots$, 397) by a function  
\begin{equation}
f_{1}(x)\doteq\frac{a_{1}}{\left(x+c_{1}\right)^{e_{1}}}+b_{1},\label{eq:C1(n)-Style}
\end{equation}
  Let $C_{1}(n)$ = $\frac{3}{2\pi^{2}}\left(\ln\left(\frac{12\sqrt{2n^{3}}h(n)}{\pi}\right)\right)^{2}-n$.
 The reason that we fit $C_{1}(x)$ by the function in the form displayed
in \eqref{eq:C1(n)-Style} is the same as that described in section
3 of \cite{liwenwei2016-Estmn-pn-arXiv} (although the data differs
distinctly). Many other types of functions have been tested, but they
 can not fit these data very well.

But here it is not valid to obtain the  constants in $f_{1}(n)$
by iteration method described in reference \cite{liwenwei2016-Estmn-pn-arXiv}.

\begin{figure}
\begin{centering}
 \includegraphics[scale=0.35]{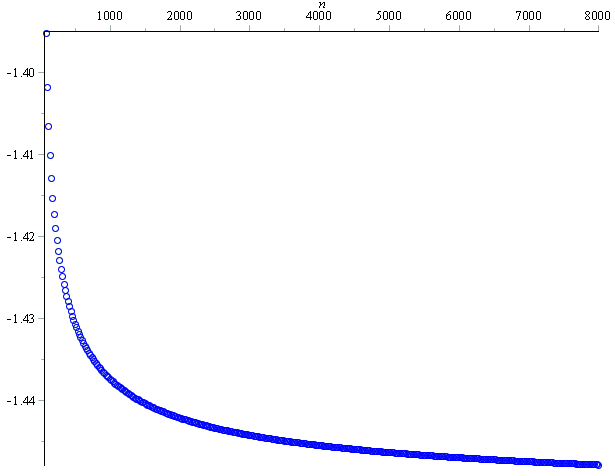}
\par\end{centering}
\caption{The graph of the data $\left(n,\frac{3}{2\pi^{2}}\left(\ln\left(\frac{12\sqrt{2n^{3}}h(n)}{\pi}\right)\right)^{2}-n\right)$}
 \label{Fig:n_C1(n)}
\end{figure}

The figure of the data $\left(n,\,\frac{3}{2\pi^{2}}\left(\ln\left(\frac{12\sqrt{2n^{3}}h(n)}{\pi}\right)\right)^{2}-n\right)$
( $n$ = $60+20k$, $k$ = 1, 2, $\cdots$, 397) is shown on Figure
\ref{Fig:n_C1(n)} on page \pageref{Fig:n_C1(n)}.

First, we try to fit $\left(n,\,\frac{3}{2\pi^{2}}\left(\ln\left(\frac{12\sqrt{2n^{3}}h(n)}{\pi}\right)\right)^{2}-n\right)$
( $n$ = $60+20k$, $k$ = 1, 2, $\cdots$, 397) by a function in
the form 
\begin{equation}
f_{1}(x)\doteq\frac{a_{1}}{\sqrt{x+c_{1}}}+b_{1}.\label{eq:C1(n)-Style-1}
\end{equation}
 That means we have assumed that $e_{1}=1/2$, temporarily. We will
explain the reason in subsection \ref{subsec:Confirm-e1}.

The average  error of $f_{1}(x)$ is 
\begin{align}
\nonumber \\
E_{1}= & \sqrt{\frac{1}{K_{1}}\sum\limits _{n}\left(C_{1}(n)-f_{1}(n)\right)^{2}}\\
= & \sqrt{\frac{1}{K_{1}}\sum\limits _{n}\left(C_{1}(n)-\frac{a_{1}}{\sqrt{n+c_{1}}}-b_{1}\right)^{2}}\nonumber 
\end{align}
 $E_{1}=\sqrt{\frac{1}{K_{1}}\sum\limits _{k=1}^{K_{1}}\left(C_{1}(60+20k)-\frac{a_{1}}{\sqrt{60+20k+c_{1}}}-b_{1}\right)^{2}}$,\\
 where $K_{1}$ = 397, $n$ ranges from 80 to 8000, by step 20. Here
only $a_{1}$, $b_{1}$, and $c_{1}$ are unknown, so we can consider
$E_{1}$ as a function of the variables $\left(a_{1},b_{1},c_{1}\right)$.

We want to find a triple $\left(a_{1},b_{1},c_{1}\right)$ such that
$E_{1}$ reaches its minimum, or to make $E_{1}$ as small as possible. 

Since a lot of functions have several local minimum points, it is
necessary to find out whether $E_{1}$ = $E_{1}\left(a_{1},b_{1},c_{1}\right)$
has more than one local minimum before we start to calculate the minimum
point by numeral method. But $E_{1}$ = $E_{1}\left(a_{1},b_{1},c_{1}\right)$
is too complicate, it is very difficult to know all the critical points
in the range we are considering. 

\subsubsection{Preparation work\label{subsec:Preparation-work}}

For a given pair $\left(a_{1},c_{1}\right)$, by the property of the
arithmetic mean, \footnote{$\ $ For some given data $x_{i}$ ($i$ = 1, 2, $\cdots$, $k$;
$x_{i}\in\mathbb{R}$), the function $s(t)=\sqrt{\frac{1}{k}\sum\limits _{i=1}^{k}\left(x_{i}-t\right)^{2}}$
reaches its minimum at $t$ = $\frac{1}{k}\sum\limits _{i=1}^{k}x_{i}$.} it is clear that $E_{1}$ reaches its minimum when 
\begin{align}
b_{1} & =\frac{1}{K_{1}}\sum\limits _{n}\left(C_{1}(n)-\frac{a_{1}}{\sqrt{n+c_{1}}}\right)\nonumber \\
 & =\overline{C_{1}}-\frac{1}{K_{1}}\sum\limits _{n'}\frac{a_{1}}{\sqrt{n'+c_{1}}},
\end{align}
   where $\overline{C_{1}}$ = $\frac{1}{K_{1}}\sum\limits _{n'}C_{1}(n')$.

Let \\
$E_{1}=\sqrt{\frac{1}{K_{1}}\sum\limits _{n}\left(C_{1}(n)-\frac{a_{1}}{\sqrt{n+c_{1}}}-\overline{C_{1}}+\frac{1}{K_{1}}\underset{_{n'}}{\sum}\frac{a_{1}}{\sqrt{n'+c_{1}}}\right)^{2}}$
 be the average error of the the fitting function $f_{1}(x)\doteq\frac{a_{1}}{\sqrt{x+c_{1}}}+\overline{C_{1}}-\frac{1}{K_{1}}\sum\limits _{n'}\frac{a_{1}}{\sqrt{n'+c_{1}}}$.
(Here $a_{1}$ and $b_{1}$ are undetermined coefficients.)

\begin{figure*}[!tp]
\centering{}%
\begin{minipage}[t]{0.45\textwidth}%
\begin{center}
\centering \includegraphics[scale=0.34]{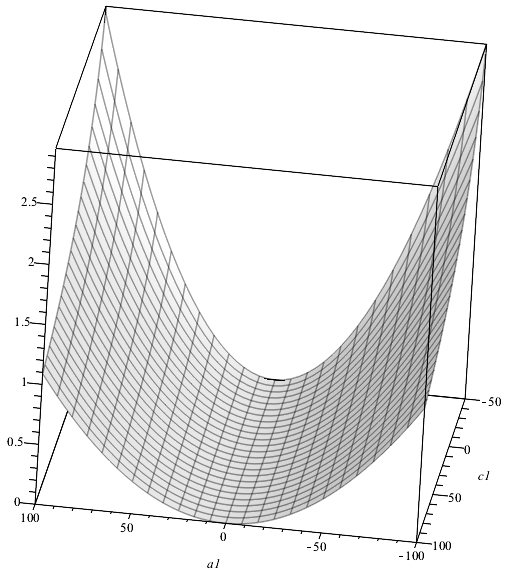}
\par\end{center}
\begin{center}
\caption{The graph of the data $\left(a_{1},c_{1},G_{1}\right)$ when $-100$
$\leqslant$ $a_{1}$ $\leqslant$ 100, $-50$ $\leqslant$ $c_{1}$
$\leqslant$ 100 }
 \label{Fig:a1_c1_G1_One}
\par\end{center}%
\end{minipage}\qquad{}%
\begin{minipage}[t]{0.45\textwidth}%
\begin{center}
\centering\includegraphics[scale=0.28]{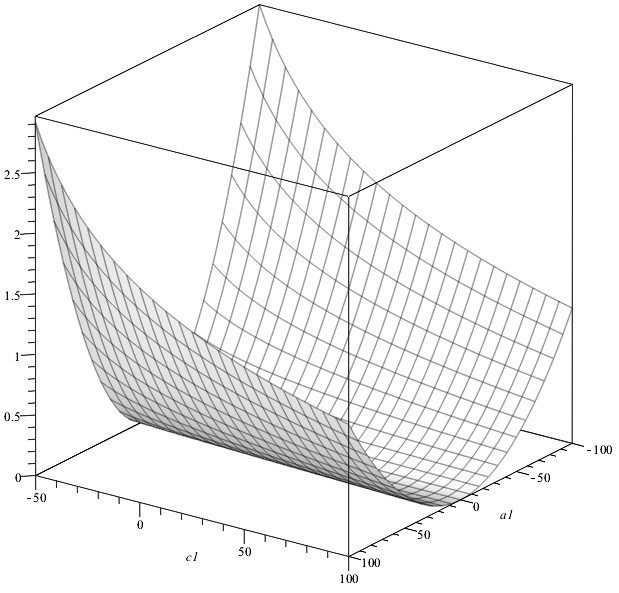}
\par\end{center}
\begin{center}
\caption{The graph of the data $\left(a_{1},c_{1},G_{1}\right)$ when $-100$
$\leqslant$ $a_{1}$ $\leqslant$ 100, $-50$ $\leqslant$ $c_{1}$
$\leqslant$ 100 version 2}
 \label{Fig:a1_c1_G1_One-B}
\par\end{center}%
\end{minipage}\vspace{0.5cm}
\begin{minipage}[t]{0.45\textwidth}%
\centering\includegraphics[scale=0.39]{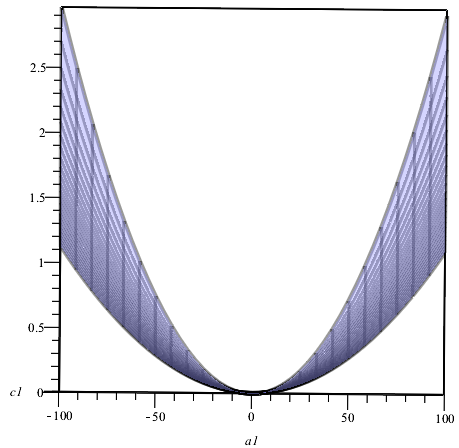}
\begin{center}
\caption{The projection of the graph of the data $\left(a_{1},c_{1},G_{1}\right)$
on the $a_{1}-G_{1}$ plane when $-100$ $\leqslant$ $a_{1}$ $\leqslant$
100, $-50$ $\leqslant$ $c_{1}$ $\leqslant$ 100}
 \label{Fig:a1_c1_G1_One-C}
\par\end{center}%
\end{minipage}\qquad{}%
\begin{minipage}[t]{0.45\textwidth}%
\centering\includegraphics[scale=0.39]{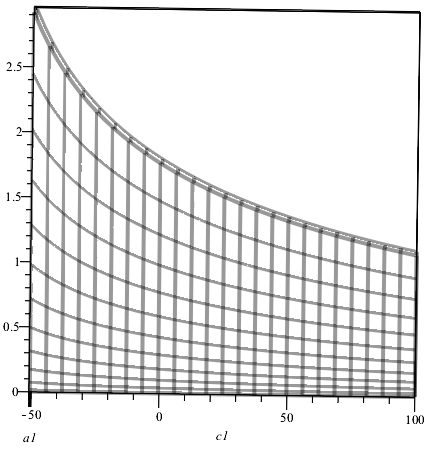}
\begin{center}
\caption{The projection of the graph of the data $\left(a_{1},c_{1},G_{1}\right)$
on the $c_{1}-G_{1}$ plane when $-100$ $\leqslant$ $a_{1}$ $\leqslant$
100, $-50$ $\leqslant$ $c_{1}$ $\leqslant$ 100}
 \label{Fig:a1_c1_G1_One-D}
\par\end{center}%
\end{minipage}
\end{figure*}

\begin{figure*}[!tp]
\begin{minipage}[t]{0.45\textwidth}%
\centering\includegraphics[scale=0.3]{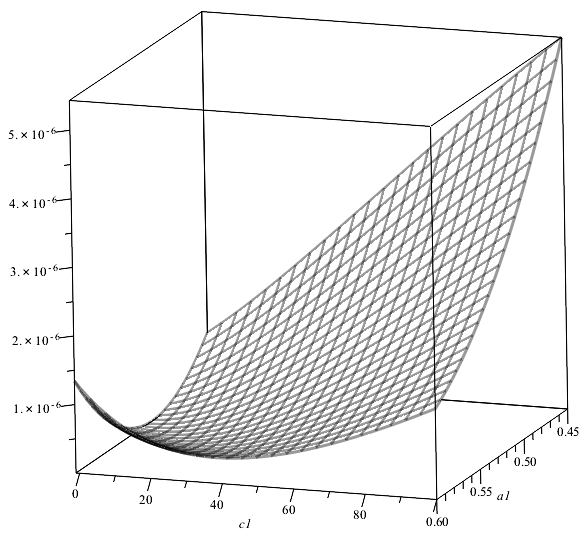}
\begin{center}
\caption{The graph of the data $\left(a_{1},c_{1},G_{1}\right)$ when $0.45$
$\leqslant$ $a_{1}$ $\leqslant$ 0.60, $-1$ $\leqslant$ $c_{1}$
$\leqslant$ 100}
 \label{Fig:a1_c1_G1_Six-A}
\par\end{center}%
\end{minipage}\qquad{}%
\begin{minipage}[t]{0.45\textwidth}%
\centering\includegraphics[scale=0.3]{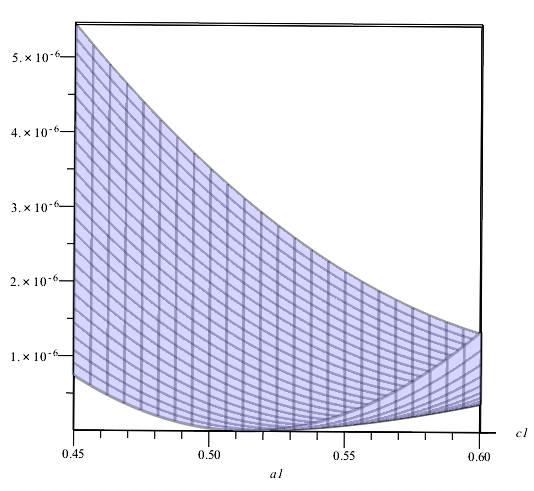}
\begin{center}
\caption{The projection of the graph of the data $\left(a_{1},c_{1},G_{1}\right)$
on the $a_{1}-G_{1}$ plane when $0.45$ $\leqslant$ $a_{1}$ $\leqslant$
0.60, $-1$ $\leqslant$ $c_{1}$ $\leqslant$ 100}
 \label{Fig:a1_c1_G1_Six-B}
\par\end{center}%
\end{minipage}\vspace{0.5cm}

\begin{minipage}[t]{0.45\textwidth}%
\centering\includegraphics[scale=0.33]{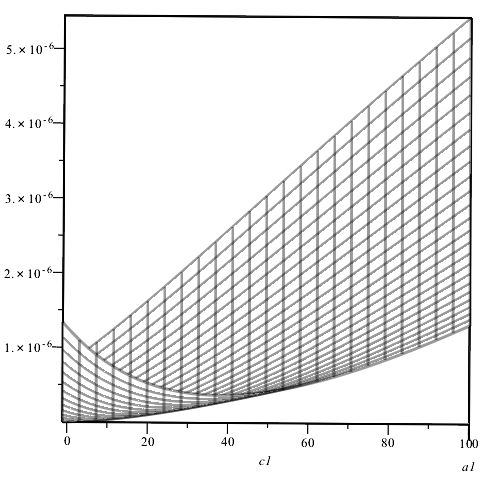}
\begin{center}
\caption{The projection of the graph of the data $\left(a_{1},c_{1},G_{1}\right)$
on the $c_{1}-G_{1}$ plane when $0.45$ $\leqslant$ $a_{1}$ $\leqslant$
0.60, $-1$ $\leqslant$ $c_{1}$ $\leqslant$ 100}
 \label{Fig:a1_c1_G1_Six-C}
\par\end{center}%
\end{minipage}\qquad{}%
\begin{minipage}[t]{0.45\textwidth}%
\centering\includegraphics[scale=0.46]{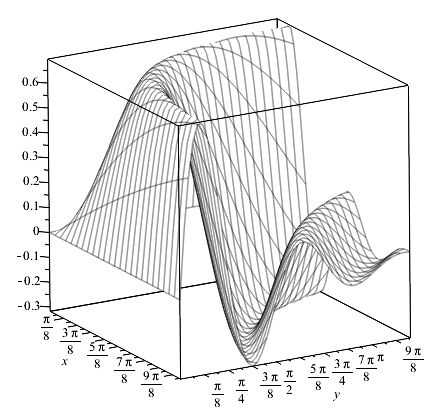}
\begin{center}
\caption{Example of a surface with more than one bottom. }
 \label{Fig:Examp-surface-more-bottom}
\par\end{center}%
\end{minipage}
\end{figure*}

Let $G_{1}=E_{1}^{2}=$\\
 $\frac{1}{K_{1}}\sum\limits _{n}\left(C_{1}(n)-\frac{a_{1}}{\sqrt{n+c_{1}}}-\overline{C_{1}}+\frac{1}{K_{1}}\sum\limits _{n'}\frac{a_{1}}{\sqrt{n'+c_{1}}}\right)^{2}.$
 Here only $a_{1}$ and $b_{1}$ are unknowns. $G_{1}$ could be
considered as a function of $\left(a_{1},c_{1}\right)$. In order
to find the minimum point of $G_{1}$,  we can draw the figure of
the function $G_{1}=G_{1}\left(a_{1},c_{1}\right)$, (In a cube coordinate
system with axis $a_{1},c_{1}$ and $G_{1}$) as shown on Figure \ref{Fig:a1_c1_G1_One},
Figure \ref{Fig:a1_c1_G1_One-B} and Figure \ref{Fig:a1_c1_G1_Six-A}.
 Figure \ref{Fig:a1_c1_G1_One-C} and Figure \ref{Fig:a1_c1_G1_One-D}
are the projection of the graph of $\left(a_{1},c_{1},G_{1}\right)$
(when $-100$ $\leqslant$ $a_{1}$ $\leqslant$ 100, $-50$ $\leqslant$
$c_{1}$ $\leqslant$ 100, which is a part of a surface) on the $a_{1}-G_{1}$
plane (spanned by the axis $a_{1}$ and $G_{1}$) and $c_{1}-G$ plane,
respectively.

From these figures, we can find out that the influence of $c_{1}$
to $G_{1}$ is much less than that of $a_{1}$. In Figure \ref{Fig:a1_c1_G1_Six-C},
we find that when $G_{1}$ reaches its minimum, $a_{1}$ is between
0.50 and 0.53, but there is not a definite range for $c_{1}$. 

It is possible that for different range of $c_{1}$, the range of
$a_{1}$ when $G_{1}$ reaches its minimum will be different. But
considering that $\frac{a_{1}}{\sqrt{n+c_{1}}}+b_{1}$ is a real,
$c_{1}$ should be greater than $-1$ in theory. (For the fitting
data used here, $c_{1}$ should be greater than $-80$.) From Figure
\ref{Fig:a1_c1_G1_One-C}, we can see that $G_{1}$ touches its bottom
when $-15$ $\leqslant$ $a_{1}$ $\leqslant$ 15. Although we can
not see clearly the exact value of of $a_{1}$ in the minimum points,
we can draw another figure of $\left(a_{1},c_{1},G_{1}\right)$ when
$-15$ $\leqslant$ $a_{1}$ $\leqslant$ 15 and $-1$ $\leqslant$
$c_{1}$ $\leqslant$ 100 to observe more details (the figure is not
presented here), then we will find that the more detailed range of
$a_{1}$ for the minimum points is $[-3,3]$ in the new figure (not
presented), next, draw the figures of $\left(a_{1},c_{1},G_{1}\right)$
when $-3$ $\leqslant$ $a_{1}$ $\leqslant$ 3, $0$ $\leqslant$
$a_{1}$ $\leqslant$ 1, $0.2$ $\leqslant$ $a_{1}$ $\leqslant$
0.8 or $0.45$ $\leqslant$ $a_{1}$ $\leqslant$ 0.6, respectively,
while $-1$ $\leqslant$ $c_{1}$ $\leqslant$ 100, we will find the
range of $a_{1}$ of the minimum points of $G_{1}$ is about {[}0.50,
0.53{]}. The projections of the figure \ref{Fig:a1_c1_G1_Six-A} of
$\left(a_{1},c_{1},G_{1}\right)$ when $0.45$ $\leqslant$ $a_{1}$
$\leqslant$ 0.6 and $-1$ $\leqslant$ $c_{1}$ $\leqslant$ 100
is shown on Figure \ref{Fig:a1_c1_G1_Six-B} and Figure \ref{Fig:a1_c1_G1_Six-C}.

In Figure \ref{Fig:a1_c1_G1_One-D},  for the curves on the bottom,
$G_{1}$ decrease with $c_{1}$ at first then increases with $c_{1}$,
but it is difficult to find the critical point in the figure, since
different curves have different critical points.

Although we can not find an satisfying value of $a_{1}$ or $c_{1}$
from the figures above to construct a fitting function $f_{1}$, 
these pictures show that the figure of $G_{1}=\left(a_{1},c_{1},G_{1}\right)$
has only one bottom in the domain we are considering, unlike the figure
of another function shown on Figure \ref{Fig:Examp-surface-more-bottom},
so the existence of the minimum point is almost confirmed, therefore
we are confident to find the value of $a_{1}$ or $c_{1}$ in the
minimum point by numerical method.  This guarantees the validity
of the numerical calculation by loop in next step. 

\subsubsection{Find $c_{1}$  \label{subsec:Find-c1}}

On the other hand, by the least square method, to fit the data $\left(x_{k},\,y_{k}\right)$
($k$ = 1, 2, $\cdots$, $K_{1}$) by a linear function $y=a\times x+b$,
the result is that 
\begin{eqnarray}
a & = & \frac{\overline{x_{n}\cdot y_{n}}-\overline{x_{n}}\cdot\overline{y_{n}}}{\overline{x_{n}^{2}}-\left(\overline{x_{n}}\right)^{2}},\label{eq:coeff-a}\\
b & = & \frac{\overline{x_{n}}\cdot\overline{x_{n}\cdot y_{n}}-\overline{x_{n}^{2}}\cdot\overline{y_{n}}}{\left(\overline{x_{n}}\right)^{2}-\overline{x_{n}^{2}}}=\overline{y_{n}}-\overline{x_{n}}\cdot a,\label{eq:coeff-b}
\end{eqnarray}
 where $\overline{x_{n}}$ = $\frac{1}{K_{1}}\sum\limits _{i=1}^{K_{1}}x_{i}$
is the average value of $x_{i}$,  $\overline{y_{n}}$ = $\frac{1}{K_{1}}\sum\limits _{i=1}^{K_{1}}y_{i}$,
$\overline{x_{i}^{2}}$ = $\frac{1}{K_{1}}\sum\limits _{i=1}^{K_{1}}x_{i}^{2}$,
$\overline{x_{n}\cdot y_{n}}$ = $\frac{1}{K_{1}}\sum\limits _{i=1}^{K_{1}}x_{i}\cdot y_{i}$.
By definition, $\left(\overline{x_{n}}\right)^{2}$ = $\left(\frac{1}{K_{1}}\sum\limits _{i=1}^{K_{1}}x_{i}\right)^{2}$
is the square of the average value of $x_{n}$. So, by the least square
method, $a$ and $b$ are uniquely determined by the given data $\left(x_{k},\,y_{k}\right)$
($k$ = 1, 2, $\cdots$, $K_{1}$).

For every given value of $c_{1}$ (greater than $-80$), we can fit
$\left(n,\,C_{1}(n)\right)$ ($n$ = $60+20k$, $k$ = 1, 2, $\cdots$,
397) by a function $f_{1}(x)\doteq\frac{a_{1}}{\sqrt{x+c_{1}}}+b_{1}$
by the least square method, if we consider $\frac{1}{\sqrt{60+20k+c_{1}}}$
and $C_{1}(60+20k)$ as $x_{k}$ and $y_{k}$, respectively.  Then 

$\overline{x_{n}}$ = $\frac{1}{K_{1}}\sum\limits _{k=1}^{K_{1}}x_{k}$
= $\frac{1}{K_{1}}\sum\limits _{k=1}^{K_{1}}\frac{1}{\sqrt{60+20k+c_{1}}}$,
$\negthickspace$ $\overline{y_{n}}$ = $\frac{1}{K_{1}}\sum\limits _{k=1}^{K_{1}}y_{k}$
= $\frac{1}{K_{1}}\sum\limits _{k=1}^{K_{1}}C_{1}(60+20k)$,

$\overline{x_{n}^{2}}$ = $\frac{1}{K_{1}}\sum\limits _{k=1}^{K_{1}}x_{k}^{2}$
= $\frac{1}{K_{1}}\sum\limits _{k=1}^{K_{1}}\frac{1}{60+20k+c_{1}}$,

$\overline{x_{n}\cdot y_{n}}$ = $\frac{1}{K_{1}}\sum\limits _{k=1}^{K_{1}}x_{k}\cdot y_{k}$
= $\frac{1}{K_{1}}\sum\limits _{k=1}^{K_{1}}\frac{C_{1}(60+20k)}{\sqrt{60+20k+c_{1}}}$,

$a_{1}=\frac{\overline{x_{n}\cdot y_{n}}-\overline{x_{n}}\cdot\overline{y_{n}}}{\overline{x_{n}^{2}}-\left(\overline{x_{n}}\right)^{2}}$,
$\quad$ $b_{1}=\overline{y_{n}}-\overline{x_{n}}\cdot a_{1}$.

So $a_{1}$ and $b_{1}$  could both be considered as functions
of $c_{1}$, denoted by $a_{1}=a_{1}(c_{1})$, $b_{1}=b_{1}(c_{1})$,
since they are uniquely determined by $c_{1}$ with the given data.

Then $G_{2}=E_{1}^{2}=$ \\
$\frac{1}{K_{1}}\sum\limits _{k=1}^{K_{1}}\left(C_{1}(60+20k)-\frac{a_{1}(c_{1})}{\sqrt{60+20k+c_{1}}}-b_{1}(c_{1})\right)^{2}$
is a  function of $c_{1}$.

It  will cost some time to plot the figure of the function $G_{2}$
= $G_{2}(c_{1})$ in a CAS software.

\begin{figure}
\centering\includegraphics[scale=0.27]{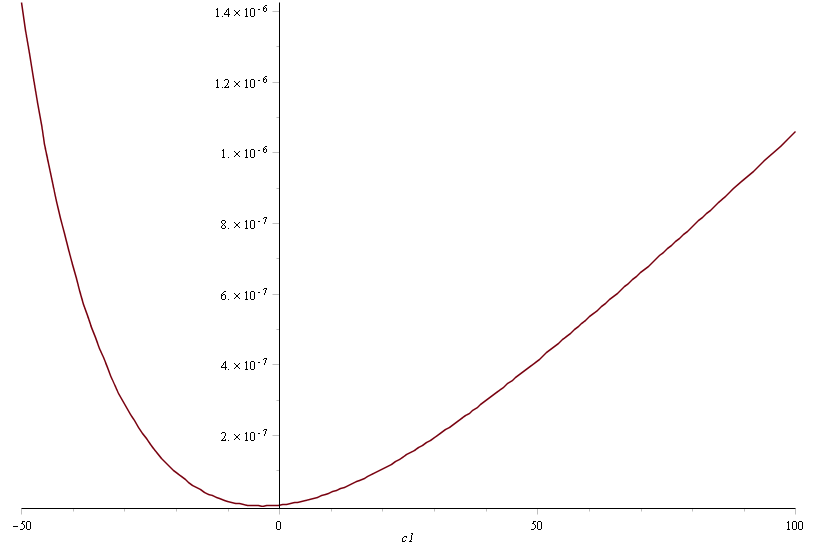}
\begin{centering}
\caption{The graph of the function $G_{2}$ = $G_{2}(c_{1})$ when $-50$ $\leqslant$
$c_{1}$ $\leqslant$ 100}
 \label{Fig:c1_G2_A}\vspace{0.5cm}
\par\end{centering}
\centering\includegraphics[scale=0.28]{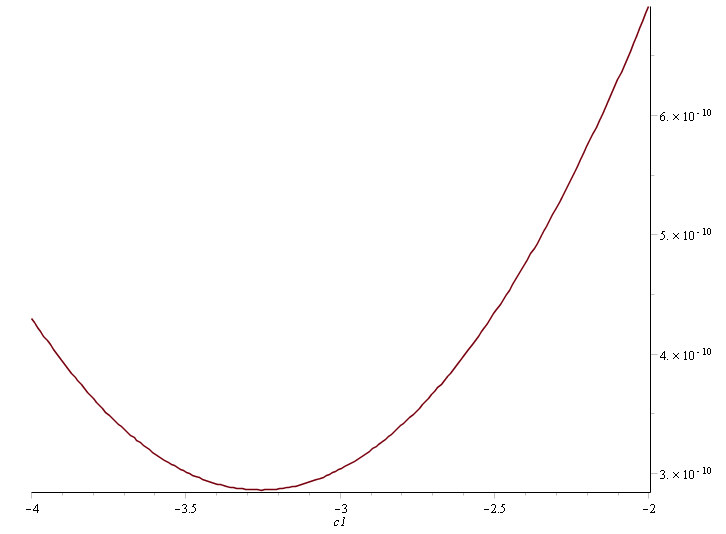}
\begin{centering}
\caption{The graph of the function $G_{2}$ = $G_{2}(c_{1})$ when $-4$ $\leqslant$
$c_{1}$ $\leqslant$ -2}
 \label{Fig:c1_G2_B}
\par\end{centering}
\begin{centering}
\vspace{0.5cm}
\par\end{centering}
\begin{centering}
\includegraphics[scale=0.21]{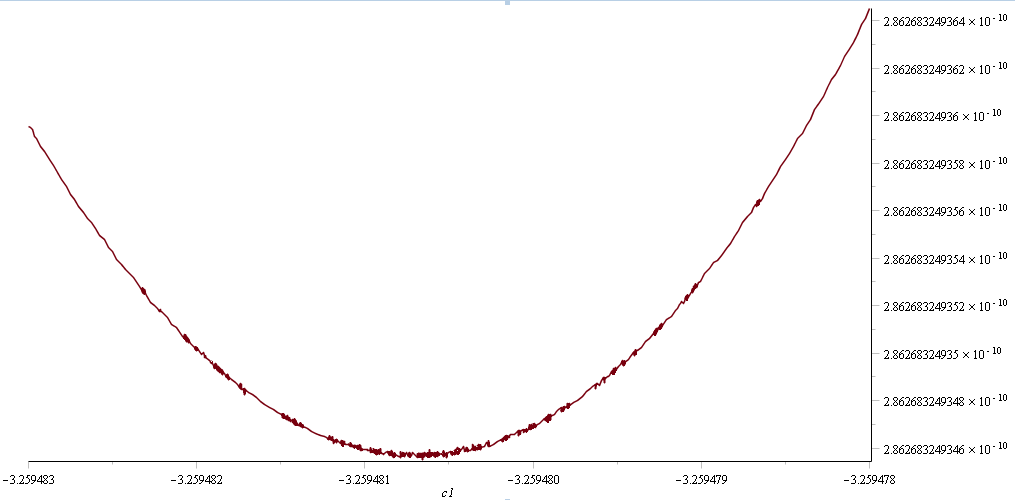}
\par\end{centering}
\centering{}\caption{The graph of the function $G_{2}$ = $G_{2}(c_{1})$ when $-3.259483$
$\leqslant$ $c_{1}$ $\leqslant$ $-3.259478$}
 \label{Fig:c1_G2_C}
\end{figure}

If we plot the figure of the function $G_{2}$ = $G_{2}(c_{1})$ on
the coordinates (as shown on Figure \ref{Fig:c1_G2_A}, Figure \ref{Fig:c1_G2_B}
and Figure \ref{Fig:c1_G2_C}), we will find that $G_{2}$ reaches
its minimum when $c_{1}$ $\approx$ $-3.2594807$. On Figure \ref{Fig:c1_G2_C},
we find that the curve of $G_{2}$ = $G_{2}(c_{1})$ is not so smooth.
The reason is that we hold up 18 significant digits in the process.
If we compute more significant digits in the process, the curve on
Figure \ref{Fig:c1_G2_C} will be more smooth, at the cost of much
more time. By writing a small program (since the default function
to find the minimum provide by the software Maple 18 are unable to
deal with such a complicated function $G_{2}$ = $G_{2}(c_{1})$ involving
so much data), we can obtain a more accurate value of the critical
point 
\[
c_{1}=-3.259480684.
\]
  When the value of $c_{1}$ is obtained, we can find the value of
$a_{1}$ and $b_{1}$ by the least square method without difficulty,
i.e., 
\begin{align*}
a_{1} & =0.5097429624,\\
b_{1} & =-1.453552800.
\end{align*}

But here $c_{1}$ is less than $-1$,  so the estimation formula
for $h(n)$ constructed from these coefficients is invalid when $n<4$.

\subsubsection{Confirm $e_{1}$ \label{subsec:Confirm-e1}}

In \cite{liwenwei2016-Estmn-pn-arXiv} we fit $C_{1}(n)$ = $\frac{3}{2}\cdot\frac{\left(\ln\left(4n\sqrt{3}p(n)\right)\right)^{2}}{\pi^{2}}-n$
by by a function $f_{1}(x)\doteq\frac{a_{1}}{\left(n+c_{1}\right)^{e_{1}}}+b_{1}$
when estimating $p(n)$, and found that $e_{1}$ $\approx$ $0.50$
by iteration. Here the iteration method does not work well, so we
fit $C_{1}(n)$ = $\frac{3}{2\pi^{2}}\left(\ln\left(\frac{12\sqrt{2n^{3}}h(n)}{\pi}\right)\right)^{2}-n$
by a function $f_{1}(x)\doteq\frac{a_{1}}{\sqrt{n+c_{1}}}+b_{1}$
directly, which means that we have assumed that $e_{1}$ $=$ $\frac{1}{2}$.
 Here we may doubt that whether $e_{1}$ $=$ $0.5$ is the best
option for us? 

Here we use the same idea described in subsection \ref{subsec:Preparation-work}.

For every pair $(e_{1},c_{1})$, we can obtain corresponding $a_{1}$
and $b_{1}$ by the least square method, just like \eqref{eq:coeff-a}
and \eqref{eq:coeff-b}, except that here $x_{k}=\frac{1}{\left(60+20k+c_{1}\right)^{e_{1}}}$.

So the square of the average error \\
$\negthickspace$ $G_{3}=E_{1}^{2}=\frac{1}{K_{1}}\sum\limits _{n}\left(C_{1}(n)-\frac{a_{1}\left(e_{1},c_{1}\right)}{\left(n+c_{1}\right)^{e_{1}}}-b_{1}\left(e_{1},c_{1}\right)\right)^{2}$
\\
 could be considered as a function of $e_{1}$ and $c_{1}$, as both
$a_{1}=a_{1}\left(e_{1},c_{1}\right)$ and $b_{1}=b_{1}\left(e_{1},c_{1}\right)$
could be expressed by certain elementary functions of $e_{1}$ and
$c_{1}$.

\begin{figure}[H]
\begin{centering}
\par\end{centering}
\begin{centering}
\includegraphics[scale=0.36]{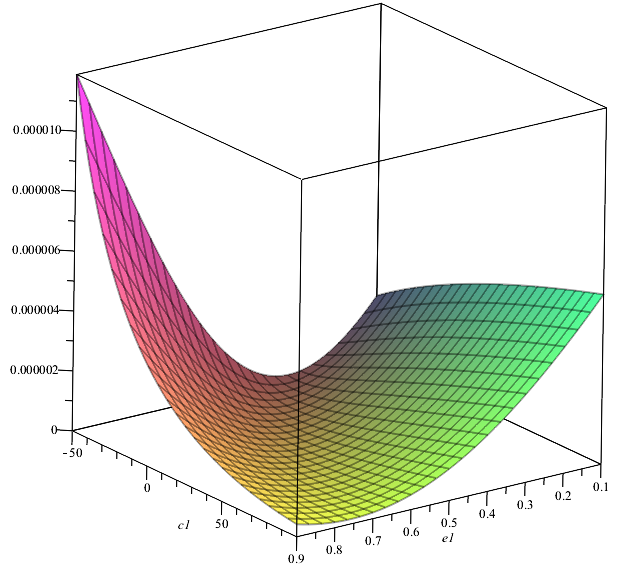}
\par\end{centering}
\caption{The graph of the function $G_{3}$ = $G_{3}(e_{1},c_{1})$ when $0.1$
$\leqslant$ $e_{1}$ $\leqslant$ $0.9$, $-50$ $\leqslant$ $c_{1}$
$\leqslant$ $100$}
 \label{Fig:e1_c1_G3}

\vspace{0.5cm}
\end{figure}

If we draw the figure of the function $G_{3}$ = $G_{3}(e_{1},c_{1})$,
we will find that the surface has only one bottom when $0.1$ $\leqslant$
$e_{1}$ $\leqslant$ $0.9$, $-50$ $\leqslant$ $c_{1}$ $\leqslant$
$100$, as shown on Figure \ref{Fig:e1_c1_G3}. But the process to
draw the figure is time-consuming. It costs more than 5 hours on a
notebook (ThinkPad E40 Edge, with 6 GB RAM and AMD P360 Dual-Core
Processor 2.30GHz) by Maple 18 in Ubuntu 14.04.1 system.

\begin{figure}[H]
\begin{centering}
\par\end{centering}
\begin{centering}
\includegraphics[scale=0.33]{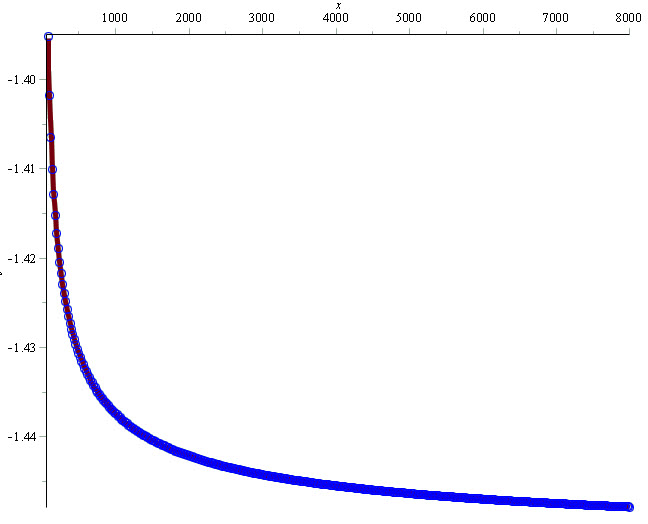}
\par\end{centering}
\caption{The graph of the data $\left(n,\frac{3}{2\pi^{2}}\left(\ln\left(\frac{12\sqrt{2n^{3}}h(n)}{\pi}\right)\right)^{2}-n\right)$
and the fitting curve }
 \label{Fig:n_C1(n)-fitting_curve}

\medskip{}

\begin{centering}
\includegraphics[scale=0.28]{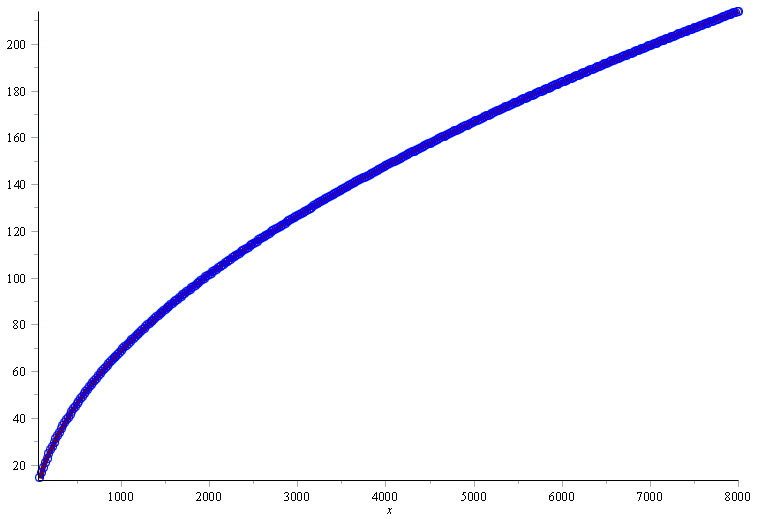}
\par\end{centering}
\caption{The graph of the data $\left(n,\ln\left(h(n)\right)\right)$ and the
fitting curve}
 \label{Fig:n_ln(h(n))-fitting_curve}

\medskip{}

\begin{centering}
\includegraphics[scale=0.33]{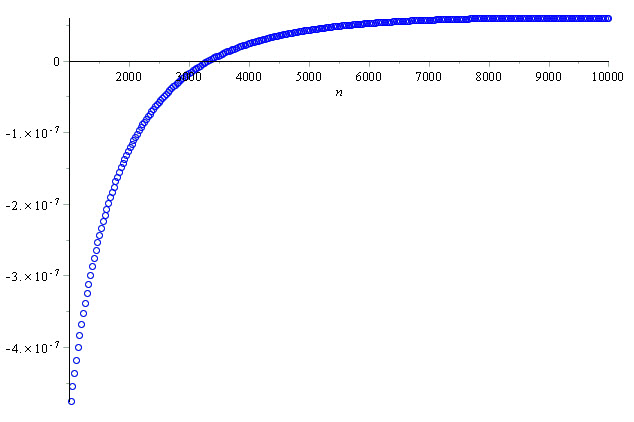}
\par\end{centering}
\caption{The Relative Error of $I_{\mathrm{ga}}(n)$when 1000 $\leqslant$
$n$ $\leqslant$ 10000, step 300}
 \label{Fig:Rel-Err-Iga}
\end{figure}

After that, by written another program, we can obtain the approximate
value of $(e_{1},c_{1})$ where $G_{3}$ touches the bottom, i.e.
$e_{1}$ $\approx$ $0.494$, $c_{1}$ $\approx$ $-4.85$, when 18
significant digits are involved in the process, which still costs
tens of minutes. Considering that we have used only a small part of
data, we can not afford the time for computing more significant digits
 in process, and the computing is so complicated hence error accumulation
effect is considerable, so we choose $e_{1}$ $=$ $0.50$ while it
differs very little with $0.494$. Another reason is that we prefer
simple exponent, as the time spend on computing a square root is much
less than that to compute a power with exponent $0.494$ in general.
 Here the value of $c_{1}\approx-4.85$ is obviously different from
the value obtained at the end of subsection \eqref{subsec:Find-c1},
because of the little difference on $e_{1}$. Therefore, it will be
fine to use the result in subsection \eqref{subsec:Find-c1}. 

\begin{figure}[H]
\begin{centering}
\includegraphics[scale=0.21]{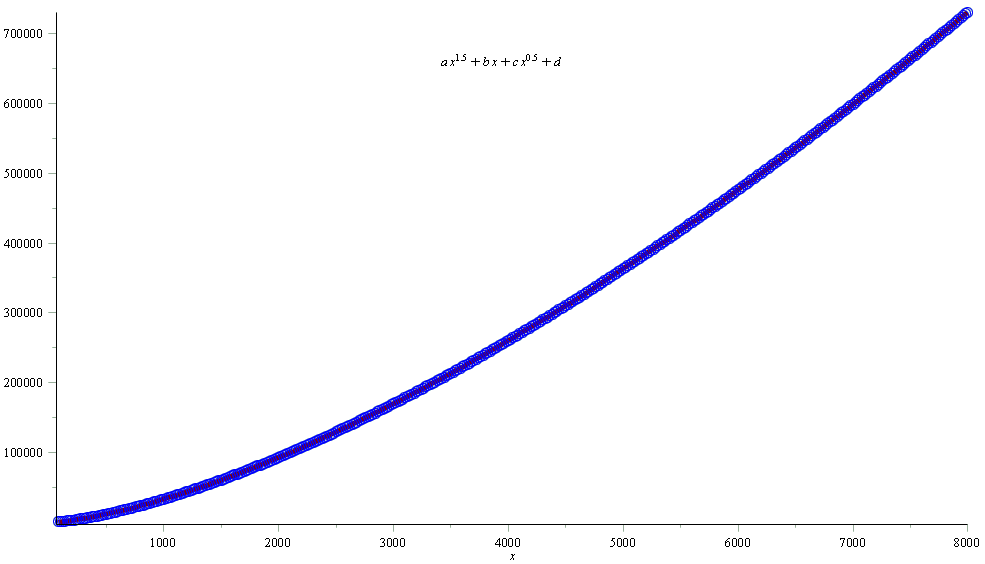}
\par\end{centering}
\begin{centering}
\caption{The graph of the data $\Bigl(n,\,\frac{\pi\exp\left(\pi\sqrt{\frac{2}{3}}\sqrt{n}\right)}{12\sqrt{2}h(n)}\Bigr)$
 and the fitting curve}
 \label{Fig:2_4_25_(n-C4(n))}
\par\end{centering}
In this figure, the points are shown as small circles which are very
close to each other. Theses crowded circles seams like a thick curve.
A fitting curve passes though the center of these circles. The fitting
curve might not be found in reduce printing. That means the curve
fits the points (displayed as circles) very well.
\end{figure}

\subsubsection{The Result}

By the value $a_{1}$ = 0.5097429624, $b_{1}$ = $-$1.453552800,
$c_{1}$ = $-3.259480684$ obtained in subsection \ref{subsec:Find-c1},
we will have a fitting function \\
$f_{1}(x)\doteq\frac{a_{1}}{\sqrt{x+c_{1}}}+b_{1}$. \\
The figure of $f_{1}(x)$ (when $-80$ $\leqslant$ $x$ $\leqslant$
8000) together with the figure of the data $\left(n,\,\frac{3}{2\pi^{2}}\left(\ln\left(\frac{12\sqrt{2n^{3}}h(n)}{\pi}\right)\right)^{2}-n\right)$
( $n$ = $60+20k$, $k$ = 1, 2, $\cdots$, 397) are shown on Figure
\ref{Fig:n_C1(n)-fitting_curve}. It shows that $f_{1}(n)$ fits $\frac{3}{2\pi^{2}}\left(\ln\left(\frac{12\sqrt{2n^{3}}h(n)}{\pi}\right)\right)^{2}-n$
very well.

Then we could fit $h(n)$ by 
\begin{equation}
I_{\mathrm{ga}}(n)=\frac{\pi\exp\left(\sqrt{\frac{2}{3}}\pi\sqrt{n+\frac{a_{1}}{\sqrt{x+c_{1}}}+b_{1}}\right)}{12\sqrt{2n^{3}}},\quad n\geqslant4.\label{eq:Iga-Expt}
\end{equation}
 The figure of the function $\ln\left(I_{\mathrm{ga}}(x)\right)$
together with the figure of the data $\left(n,\ln\left(h(n)\right)\right)$
( $n$ = $60+20k$, $k$ = 1, 2, $\cdots$, 397) are shown on Figure
\ref{Fig:n_ln(h(n))-fitting_curve}. It seems that $I_{\mathrm{ga}}(n)$
fits $h(n)$ very well. 

The relative error of $I_{\mathrm{ga}}$ is shown on Table \ref{Table:Rel-Err-h(n)-Iga(x)}
(when $n\leqslant1000$ ) and Figure \ref{Fig:Rel-Err-Iga} (when
1000 $<$ $n$ $\leqslant10000$). 

When $n<20$, the relative error of $I_{\mathrm{ga}}$ is still greater
than $2\%$. Although it is much better than the error of $I_{\mathrm{g}}$,
it is not as good as expect when $n<40$. If we take the round approximation
by 
\begin{align}
I'_{\mathrm{ga}}(n) & =\left\lfloor \frac{\pi\exp\left(\sqrt{\frac{2}{3}}\pi\sqrt{n+\frac{a_{1}}{\sqrt{x+c_{1}}}+b_{1}}\right)}{12\sqrt{2n^{3}}}+\frac{1}{2}\right\rfloor ,\nonumber \\
 & \quad n\geqslant4,\label{eq:Iga-round}
\end{align}
 the relative error will be obviously smaller with a few exceptions,
as shown on Table \ref{Table:Rel-Err-h(n)-Iga(n)-round}.

Later we will find out that it is obviously greater than the relative
error of $I_{\mathrm{g1}}$ and $I_{\mathrm{g2}}$ obtained in the
next subsection by modifying the denominator part; when $4000<n<10000$,
the relative error of $I_{\mathrm{ga}}$ is about 1000 times of that
of $I_{\mathrm{g2}}$. 

\begin{table}[H]

\begin{centering}
\includegraphics[viewport=97bp 479bp 461bp 666bp,clip,scale=0.63]{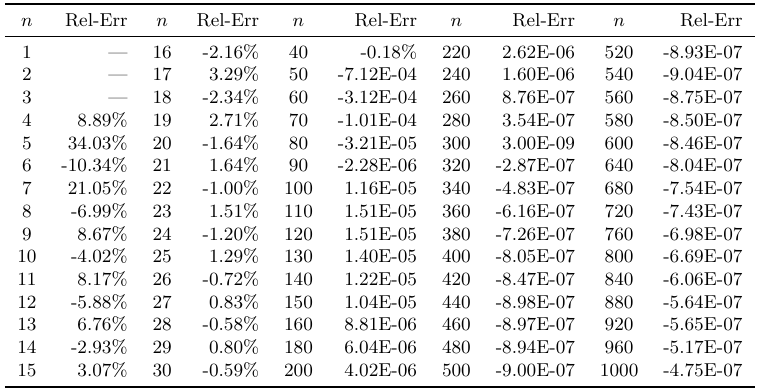}
\par\end{centering}
\caption{The relative error of $I_{\mathrm{ga}}(n)$ to $h(n)$ when $n\leqslant1000$.}
 \label{Table:Rel-Err-h(n)-Iga(x)} \vspace{0.5cm}

\begin{centering}
\includegraphics[viewport=97bp 358bp 412bp 447bp,clip,scale=0.64]{tables-v2/Table-2_4_2_0-Rel-Err-Iga_n_-Iga_n_-round}
\par\end{centering}
\caption{The relative error of $\left\lfloor I_{\mathrm{ga}}(n)+\frac{1}{2}\right\rfloor $
to $h(n)$ when $n\leqslant30$.}
 \label{Table:Rel-Err-h(n)-Iga(n)-round}When 30 $<$ $n$ $\leqslant$
1000, the relative error of $\left\lfloor I_{\mathrm{ga}}(n)+\frac{1}{2}\right\rfloor $
is closed to that of $I_{\mathrm{ga}}(n)$.

\end{table}

\subsection{Method B: Modifying the Denominator\label{subsec:Method-B}}

Since $h(n)$ $\sim$ $\frac{\pi}{12\sqrt{2n^{3}}}\exp\left(\sqrt{\frac{2}{3}}\pi\sqrt{n}\right)$,
we consider estimating $h(n)$ by \\
 $\frac{\pi}{12\sqrt{2C_{3}(n)}}\exp\left(\sqrt{\frac{2}{3}}\pi\sqrt{n}\right)$,
(i.e., fit $\frac{\pi^{2}\exp\left(2\pi\sqrt{\frac{2}{3}}\sqrt{n}\right)}{288h^{2}(n)}$
by a function $C_{3}(n)$), where $C_{3}(x)$ is a cubic function
or a function like 
\[
ax^{3}+bx^{2.5}+cx^{2}+dx^{1.5}+ex+fx^{0.5}+g.
\]
But the results are worse, as the relative errors are obviously much
greater than the relative error of $I_{\mathrm{g}}(n)$ when $n<350$.

Then we consider estimating $h(n)$ by $\frac{\pi}{12\sqrt{2}C_{4}(n)}\exp\left(\sqrt{\frac{2}{3}}\pi\sqrt{n}\right)$,
or fit $\frac{\pi\exp\left(\pi\sqrt{\frac{2}{3}}\sqrt{n}\right)}{12\sqrt{2}h(n)}$
by a function 
\begin{equation}
C_{4}(n)=a_{4}n^{1.5}+b_{4}n+c_{4}n^{0.5}+d_{4}.\label{eq:fit-h(n)-C4}
\end{equation}
The result is very good. The figure of the data $\left(n,\,\frac{\pi\exp\left(\pi\sqrt{\frac{2}{3}}\sqrt{n}\right)}{12\sqrt{2}h(n)}\right)$
and the fitting curve $C_{4}(n)$ are shown on  Figure \ref{Fig:2_4_25_(n-C4(n))}
on page \pageref{Fig:2_4_25_(n-C4(n))}. Here the fitting curve is
displayed by a thick continuous curve, which lies in the middle of
the area the circles occupied. Since the circles are too crowded,
the circles themselves look like a very thick curve.

\begin{table}[H]

\begin{centering}
\includegraphics[viewport=97bp 479bp 458bp 666bp,clip,scale=0.63]{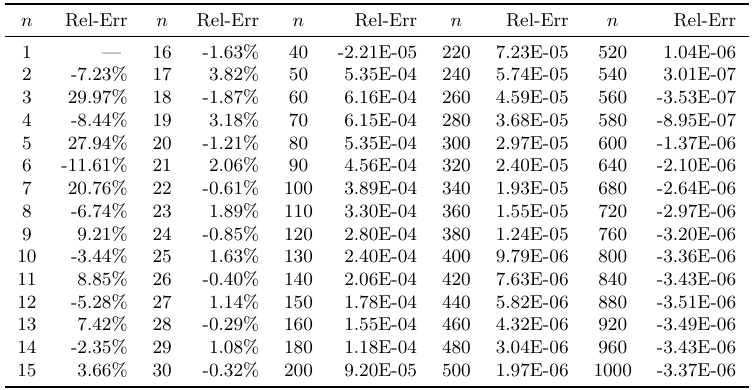}
\par\end{centering}
\caption{The relative error of $I_{\mathrm{g1}}(n)$ to $h(n)$ when $n\leqslant1000$.}
 \label{Table:Rel-Err-h(n)-Ig1(x)} \vspace{0.8cm}

\begin{centering}
\includegraphics[viewport=97bp 358bp 412bp 447bp,clip,scale=0.63]{tables-v2/Table-2_4_2_b-Rel-Err-Ig1_n_-Ig1_n_round}
\par\end{centering}
\caption{The relative error of $\left\lfloor I_{\mathrm{g1}}(n)+\frac{1}{2}\right\rfloor $
to $h(n)$ when $n\leqslant30$.}
 \label{Table:Rel-Err-h(n)-Ig1(n)-round}

\end{table}

\begin{figure}[H]
\begin{centering}
\includegraphics[scale=0.21]{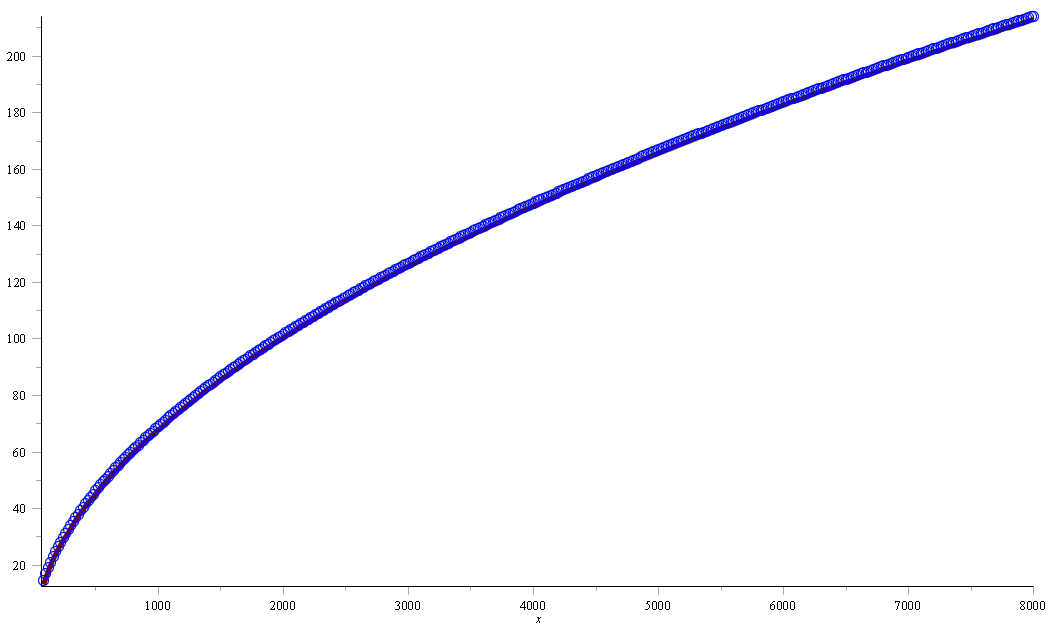}
\par\end{centering}
\caption{The graph of the data $\left(n,\ln h(n)\right)$ and the fitting curve
$\ln\left(I_{\mathrm{g}1}(n)\right)$}
 \label{Fig:2_4_26_(n-ln(h(n)))-Ig1(n)}  \vspace*{0.9cm}

\begin{centering}
\includegraphics[scale=0.21]{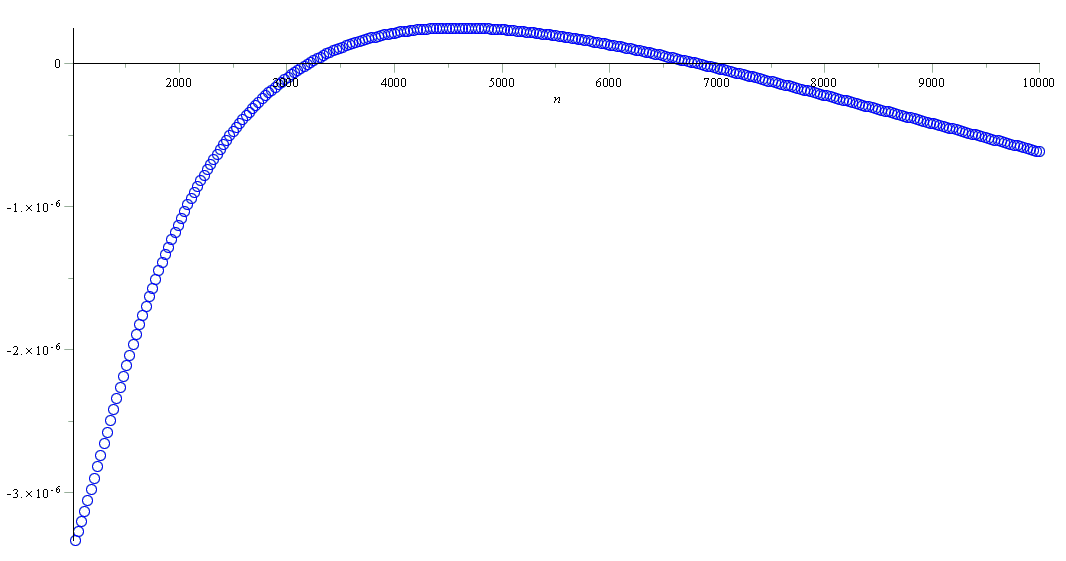}
\par\end{centering}
\caption{The Relative Error of $I_{\mathrm{g1}}(n)$when 1000 $\leqslant$
$n$ $\leqslant$ 10000, step 300}
 \label{Fig:2_4_27-Rel-Err-Ig1}

\vspace{0.6cm}

\begin{centering}
\includegraphics[scale=0.19]{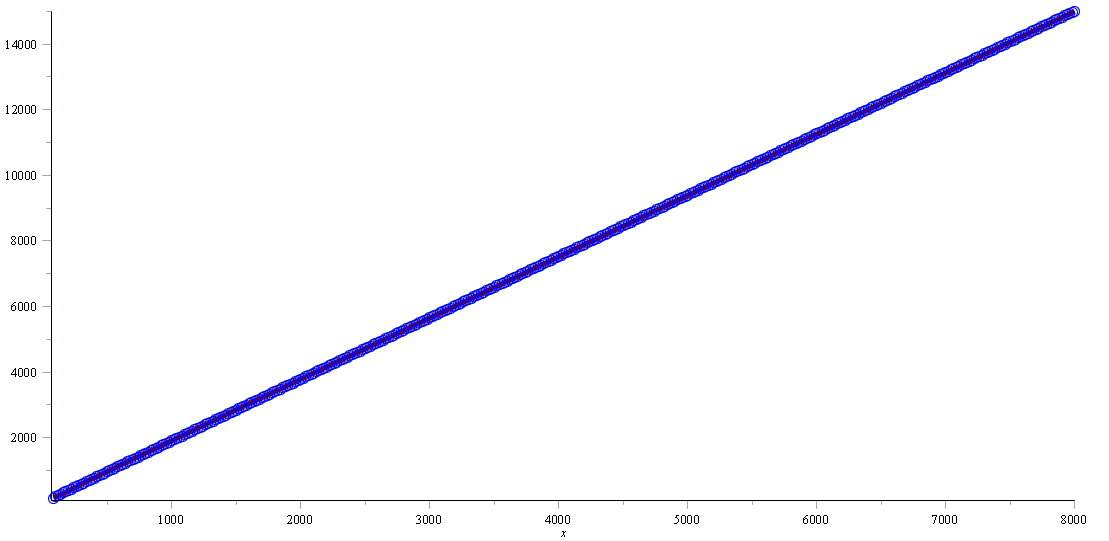}
\par\end{centering}
\begin{centering}
\caption{The graph of the data $\left(n,\ \frac{\pi\exp\left(\pi\sqrt{\frac{2}{3}}\sqrt{n}\right)}{12\sqrt{2}h(n)}-n^{3/2}\right)$
and the fitting curve $C_{5}(n)$}
 \label{Fig:2_4_28_(n,C5(n))}  \vspace{0.5cm}
\par\end{centering}
\begin{centering}
\includegraphics[scale=0.21]{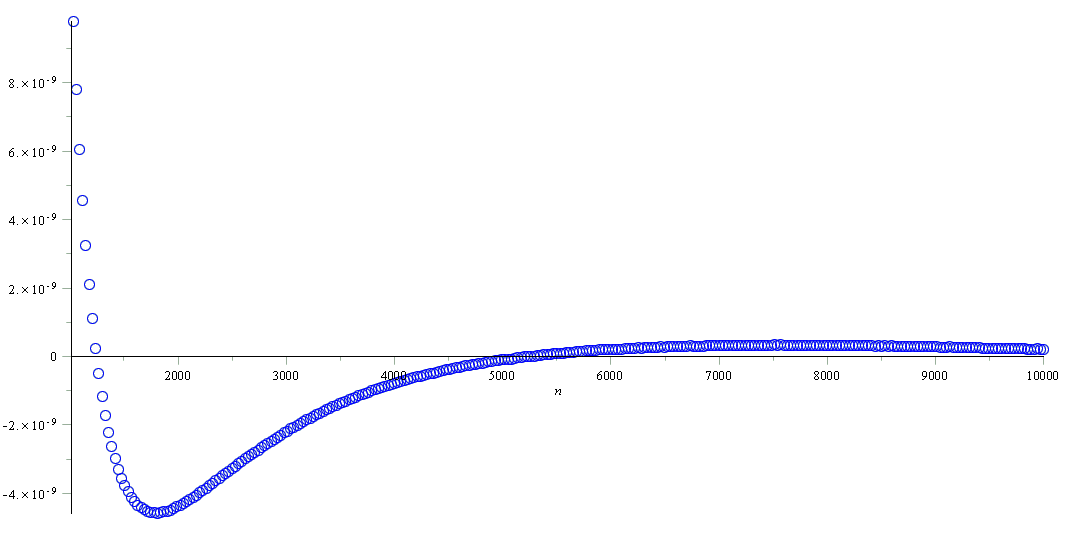}
\par\end{centering}
\centering{}\caption{The Relative Error of $I_{\mathrm{g2}}(n)$when 1000 $\leqslant$
$n$ $\leqslant$ 10000, step 300}
 \label{Fig:2_4_29-Rel-Err-Ig2}
\end{figure}

The values of the coefficients in the expression of $C_{4}(n)$ are
as follow, 
\begin{align*}
a_{4}= & 1.000010809,\\
b_{4}= & 1.862505234,\\
c_{4}= & 1.169930087,\\
d_{4}= & -0.7005460222.
\end{align*}
 The value of $a_{4}$ is very close to 1, which means that this fitting
function coincides with the Ingham-Meinardus asymptotic formula very
well.

So we have an estimation formula \label{Sym:Ig1(n)} \nomenclature[Ig1(n)]{$I_{\mathrm{g1}}(n)$}{The Ingham-Meinardus revised estimation formula 1. \pageref{Sym:Ig1(n)}}
\begin{equation}
h(n)\sim I_{\mathrm{g1}}(n)=\frac{\pi}{12\sqrt{2}C_{4}(n)}\exp\left(\sqrt{\frac{2}{3}}\pi\sqrt{n}\right).\label{eq:h(n)-estimation-1}
\end{equation}

We may call it the\emph{ Ingham-Meinardus revised estimation formula
1.} \index{Ingham-Meinardus revised estimation formula 1} The graph
of $\ln\left(I_{\mathrm{g1}}(n)\right)$ is shown on  Figure \ref{Fig:2_4_26_(n-ln(h(n)))-Ig1(n)}
on page \pageref{Fig:2_4_26_(n-ln(h(n)))-Ig1(n)}, together with the
data points of $\left(n,\ln h(n)\right)$. This revised estimation
formula is much more accurate than the asymptotic formula. The relative
error is less than $1\times10^{-6}$ when $n$ > 2000 (as shown on
Figure \ref{Fig:2_4_27-Rel-Err-Ig1} on page \pageref{Fig:2_4_27-Rel-Err-Ig1}),
 and less than $3\permil$ when $n$ $\geqslant$ 30 (as shown on
Table \ref{Table:Rel-Err-h(n)-Ig1(x)} on page \pageref{Table:Rel-Err-h(n)-Ig1(x)}).
The relative error of the round approximation $I'_{\mathrm{g1}}(n)=\left\lfloor I_{\mathrm{g1}}(n)+\frac{1}{2}\right\rfloor $
is shown on Table \ref{Table:Rel-Err-h(n)-Ig1(n)-round} on page \pageref{Table:Rel-Err-h(n)-Ig1(n)-round}.

But Equation \eqref{eq:h(n)-estimation-1} is not so satisfying when
$n<30$, especially when $n<15$ as the relative error is not negligible
for some value of $n$.

As we already know that $h(n)\sim\frac{\pi}{12\sqrt{2n^{3}}}\exp\left(\sqrt{\frac{2}{3}}\pi\sqrt{n}\right)$,
or \\
 $\ $ $n^{3/2}\sim\frac{\pi}{12\sqrt{2}h(n)}\exp\left(\sqrt{\frac{2}{3}}\pi\sqrt{n}\right)$,
which means that when fitting $\frac{\pi\exp\left(\pi\sqrt{\frac{2}{3}}\sqrt{n}\right)}{12\sqrt{2}h(n)}$
by a function $C_{4}(n)$ shown in Equation \eqref{eq:fit-h(n)-C4},
the coefficient $a_{4}$ should be exactly 1, hence we should fit
$\frac{\pi\exp\left(\pi\sqrt{\frac{2}{3}}\sqrt{n}\right)}{12\sqrt{2}h(n)}$
by a function $C'_{4}(n)=n^{3/2}+b_{5}n+c_{5}n^{1/2}+d_{5}$, or fit
$\frac{\pi\exp\left(\pi\sqrt{\frac{2}{3}}\sqrt{n}\right)}{12\sqrt{2}h(n)}-n^{3/2}$
by a function 
\begin{equation}
C_{5}(n)=b_{5}n+c_{5}n^{1/2}+d_{5}.\label{eq:fit-h(n)-C5}
\end{equation}

The figure of the data $\left(n,\,\frac{\pi\exp\left(\pi\sqrt{\frac{2}{3}}\sqrt{n}\right)}{12\sqrt{2}h(n)}-n^{3/2}\right)$
is shown on  Figure \ref{Fig:2_4_28_(n,C5(n))} on page \pageref{Fig:2_4_28_(n,C5(n))}
(together with the figure of the fitting function $C_{5}(n)$ generated
by the least square method).

The values of the coefficients in Equation \eqref{eq:fit-h(n)-C5}
are as follow 
\begin{align*}
b_{5}= & 1.864260743,\\
c_{5}= & 1.084436400,\\
d_{5}= & 0.4754177757.
\end{align*}
 So we have another estimation formula for $h(n)$, \label{Sym:Ig2(n)}
\nomenclature[Ig2(n)]{$I_{\mathrm{g2}}(n)$}{The Ingham-Meinardus revised estimation formula 2. \pageref{Sym:Ig2(n)}}
\begin{equation}
h(n)\sim I_{\mathrm{g2}}(n)=\frac{\pi\exp\left(\sqrt{\frac{2}{3}}\pi\sqrt{n}\right)}{12\sqrt{2}\left(n^{3/2}+C_{5}(n)\right)}.\label{eq:h(n)-estimation-2}
\end{equation}

\begin{table}[H]

\begin{centering}
\includegraphics[viewport=97bp 479bp 458bp 666bp,clip,scale=0.63]{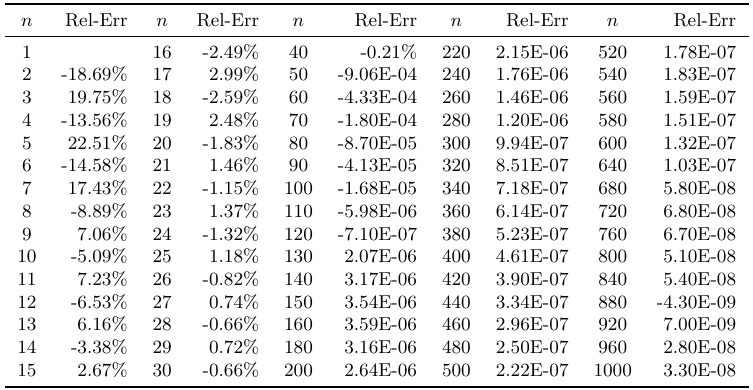}
\par\end{centering}
\caption{The relative error of $I_{\mathrm{g2}}(n)$ to $h(n)$ when $n\leqslant1000$.}
 \label{Table:Rel-Err-h(n)-Ig2(x)} \vspace{0.5cm}

\begin{centering}
\includegraphics[viewport=97bp 358bp 414bp 447bp,clip,scale=0.63]{tables-v2/Table-2_4_2_c-Rel-Err-Ig2_n_-Ig2_n_round}
\par\end{centering}
\caption{The relative error of $\left\lfloor I_{\mathrm{g2}}(n)+\frac{1}{2}\right\rfloor $
to $h(n)$ when $n\leqslant30$.}
 \label{Table:Rel-Err-h(n)-Ig2(n)-round}

\end{table}

We may call it the\emph{ Ingham-Meinardus revised estimation formula
2.} \index{Ingham-Meinardus revised estimation formula 2} The graph
of $\ln\left(I_{\mathrm{g2}}(n)\right)$ is nearly the same as that
of $\ln\left(I_{\mathrm{g1}}(n)\right)$ shown on  Figure \ref{Fig:2_4_26_(n-ln(h(n)))-Ig1(n)}
on page \pageref{Fig:2_4_26_(n-ln(h(n)))-Ig1(n)}. The second revised
estimation formula is much more accurate than the first one. The relative
error is less than $2\times10^{-9}$ when $n$ > 3000 (as shown on
Figure \ref{Fig:2_4_29-Rel-Err-Ig2}  on page \pageref{Fig:2_4_29-Rel-Err-Ig2}),
about $\frac{1}{500}$ of the relative error of $I_{\mathrm{g1}}(n)$.
When $n$ < 10, the relative error is also distinctly less than that
of $I_{\mathrm{g1}}(n)$ (as shown on Table \ref{Table:Rel-Err-h(n)-Ig2(x)}
on page \pageref{Table:Rel-Err-h(n)-Ig2(x)}). The relative error
of the round approximation $I'_{\mathrm{g2}}(n)=\left\lfloor I_{\mathrm{g2}}(n)+\frac{1}{2}\right\rfloor $
is shown on Table \ref{Table:Rel-Err-h(n)-Ig2(n)-round} (on page
\pageref{Table:Rel-Err-h(n)-Ig2(n)-round}).

It should be mentioned that in  Figure \ref{Fig:2_4_28_(n,C5(n))}
on page \pageref{Fig:2_4_28_(n,C5(n))}, the graph of the data points
lie in a line, so we might be willing to fit this line by a first
order equation. The result is 
\[
C'_{5}(n)=1.873818457\times n+27.08318017.
\]
 If we use this fitting function instead of $C_{5}(n)$ generated
above, the relative error to fit $h(n)$ will be about 10000 times
more, that is about 20 times more than that of $I_{\mathrm{g1}}(n)$.
So we did not use linear function to fit the data $\left(n,\,\frac{\pi\exp\left(\pi\sqrt{\frac{2}{3}}\sqrt{n}\right)}{12\sqrt{2}h(n)}-n^{3/2}\right)$
before. 

\begin{table}[h]

\begin{centering}
\includegraphics[viewport=97bp 515bp 456bp 703bp,clip,scale=0.63]{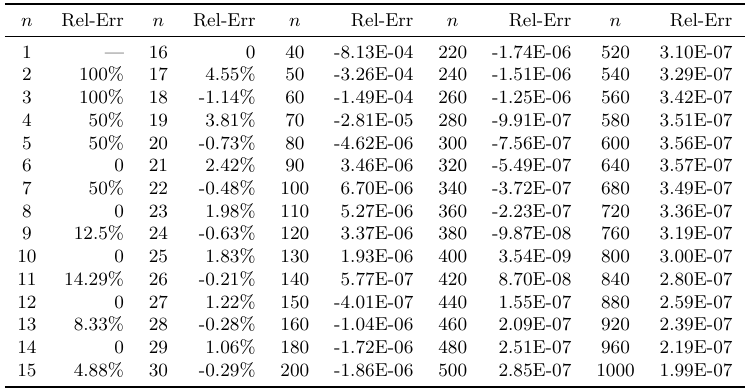}
\par\end{centering}
\caption{The relative error of $\left\lfloor F_{7a}(n)+\frac{1}{2}\right\rfloor $
to $h(n)$ when $n\leqslant1000$.}
 \label{Table:Rel-Err-h(n)-F7a(n)} \vspace{0.5cm}

\begin{centering}
\includegraphics[viewport=97bp 295bp 459bp 483bp,clip,scale=0.63]{tables-v2/Table-2_4_2_d-Rel-Err-F7-round-C7a_C7b}
\par\end{centering}
\caption{The relative error of $\left\lfloor F_{7b}(n)+\frac{1}{2}\right\rfloor $
to $h(n)$ when $n\leqslant1000$.}
 \label{Table:Rel-Err-h(n)-F7b(n)}

\end{table}

\subsection{Method C: Fittong $I_{\mathrm{g}}(n)-h(n)$\label{subsec:Method-C}}

We wonder whether we can fit $I_{\mathrm{g}}(n)-h(n)$ by a function
$r(n)$, then estimate $h(n)$ by $I_{\mathrm{g}}(n)-r(n)$ which
may be believed more accurate than $I_{\mathrm{g2}}(n)$ at the price
of being more complicated.

By the same tricks used at the beginning of this subsection, we will
have 
\[
I_{\mathrm{g}}(n)-I_{\mathrm{g}}(n-t)\sim\frac{t\pi^{2}}{24\sqrt{3}n^{2}}\exp\left(\sqrt{\frac{2}{3}}\pi\sqrt{n}\right).\quad(t\ll n)
\]

 So we may fit $I_{\mathrm{g}}(n)-h(n)$ by $\frac{\pi^{2}}{24\sqrt{3}C_{6}(n)}\exp\left(\sqrt{\frac{2}{3}}\pi\sqrt{n}\right)$
where $C_{6}(n)$ is a quadratic function or a function like 
\[
an^{2}+bn^{1.5}+cn+dn^{0.5}+e.
\]
 That means, we can fit $\frac{\pi^{2}\exp\left(\sqrt{\frac{2}{3}}\pi\sqrt{n}\right)}{24\sqrt{3}\left(I_{\mathrm{g}}(n)-h(n)\right)}$
by a function $C_{6}(n)$. But the result is useless. Although $C_{6}(n)$
will fit the data $\frac{\pi^{2}\exp\left(\sqrt{\frac{2}{3}}\pi\sqrt{n}\right)}{24\sqrt{3}\left(I_{\mathrm{g}}(n)-h(n)\right)}$
very well, but the relative error of $I_{\mathrm{g}}(n)-\frac{\pi^{2}}{24\sqrt{3}C_{6}(n)}\exp\left(\sqrt{\frac{2}{3}}\pi\sqrt{n}\right)$
to $h(n)$ is much greater than that of $I_{\mathrm{g1}}(n)$ or $I_{\mathrm{g2}}(n)$,
and the relative error differs very little with that of $I_{\mathrm{g}}(n)$
when $n$ is small. Besides, the formula $I_{\mathrm{g}}(n)-\frac{\pi^{2}}{24\sqrt{3}C_{6}(n)}\exp\left(\sqrt{\frac{2}{3}}\pi\sqrt{n}\right)$
are much more complicated than $I_{\mathrm{g1}}(n)$ and $I_{\mathrm{g2}}(n)$.

If we use the trick mentioned in \ref{subsec:Method-B}, to fit $I_{\mathrm{g}}(n)-h(n)$
by $\frac{\pi^{2}}{24\sqrt{3}\left(n^{2}+E_{6}(n)\right)}\exp\left(\sqrt{\frac{2}{3}}\pi\sqrt{n}\right)$
where $E_{6}(n)$ is a function like 
\[
bn^{1.5}+cn+dn^{0.5}+e,
\]
as we already know that the coefficient of $n^{2}$ should be 1 in
theory. The result will be a little better, but useless too. The accuracy
is not as good as that of $I_{\mathrm{g0}}(n)$.

Then we consider fitting $\frac{\pi^{2}\exp\left(\sqrt{\frac{2}{3}}\pi\sqrt{n}\right)}{24\sqrt{3}n^{2}\left(I_{\mathrm{g}}(n)-h(n)\right)}$
by a function $C_{7}(n)$. If $C_{7}(n)$ is in the form $\frac{a}{n}+b$
or $\frac{a}{n}+\frac{b}{n^{2}}+c$, the result is useless, too. If
$C_{7}(n)$ is in the form $\frac{a}{n^{0.5}}+b$, it will be barely
satisfactory. If $C_{7}(n)$ is in the form $\frac{a}{n^{0.5}}+\frac{b}{n}+\frac{c}{n^{1.5}}+\frac{d}{n^{2}}+e$
or $\frac{a}{n^{0.5}}+\frac{b}{n}+\frac{c}{n^{1.5}}+e$, the result
will be much better than the previous forms, but the accuracy (when
estimating $h(n)$) is not as good as that of $I_{\mathrm{g1}}(n)$
and $I_{\mathrm{g2}}(n)$. 

The result of $C_{7}(n)$ is 

\begin{align*}
C_{7a}(n)= & \frac{0.8782296151}{n^{0.5}}+\frac{0.2567016063}{n}\\
 & -\frac{3.580442785}{n^{1.5}}+\frac{21.28305831}{n^{2}}\\
 & +0.6879945549,
\end{align*}

or 
\begin{align*}
C_{7b}(n)= & \frac{0.8861039149}{n^{0.5}}-\frac{0.05719053203}{n}+\\
 & \frac{0.9843423289}{n^{1.5}}+0.6879343652.
\end{align*}

\begin{figure}[H]
\begin{centering}
\includegraphics[scale=0.21]{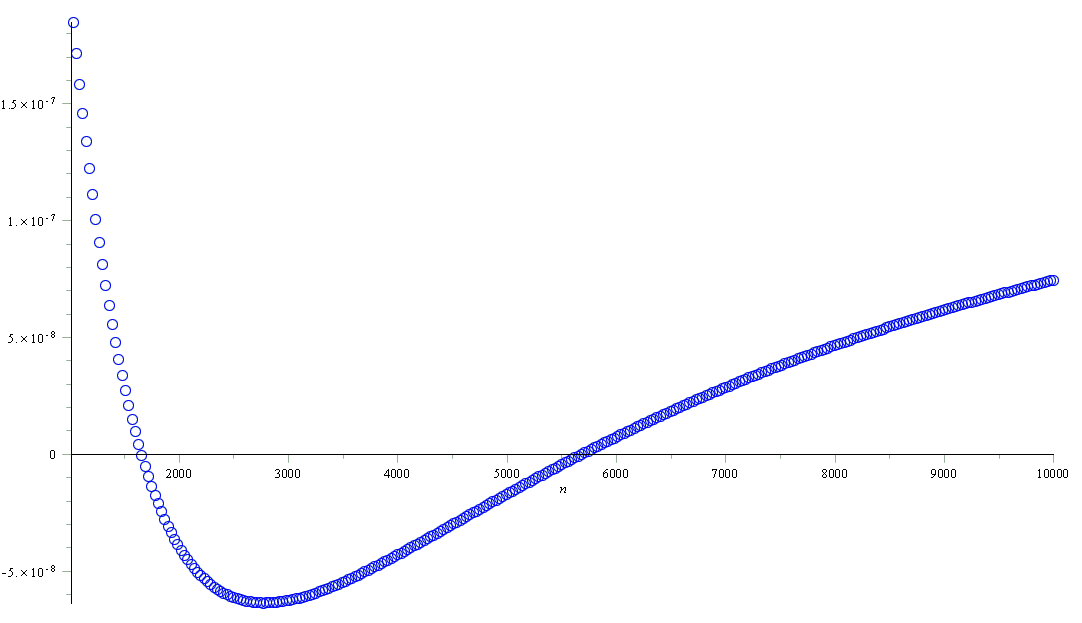}
\caption{The Relative Error of $F_{7a}(n)$ when 1000 $\leqslant$ $n$ $\leqslant$
10000, step 300}
 \label{Fig:2_4_30-Rel-Err-F7-C7a}  \vspace{0.6cm}
\par\end{centering}
\begin{centering}
\includegraphics[scale=0.2]{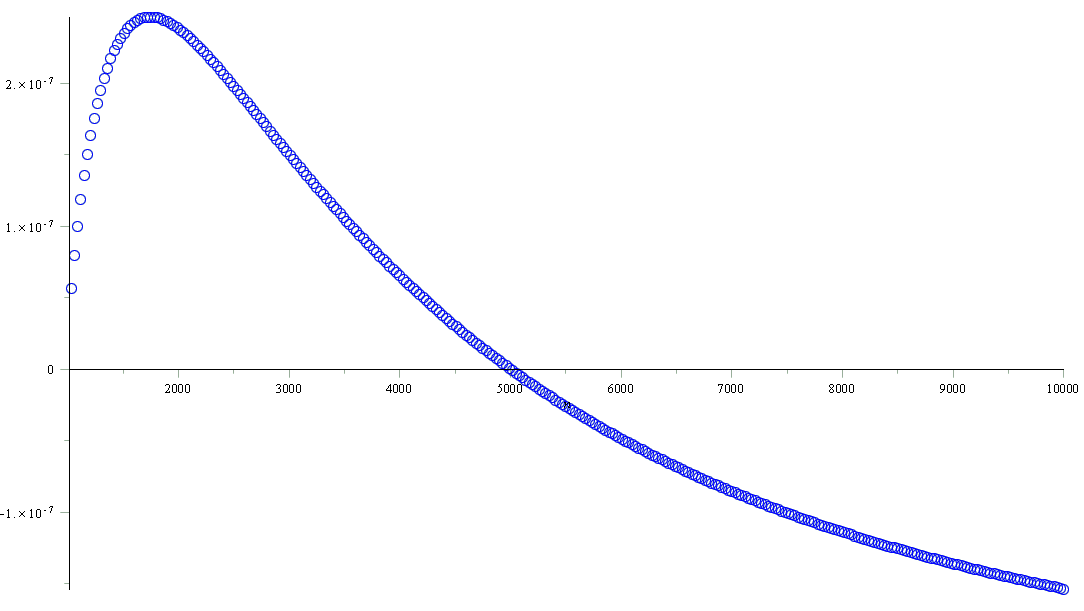}
\caption{The Relative Error of $F_{7b}(n)$ when 1000 $\leqslant$ $n$ $\leqslant$
10000, step 300}
 \label{Fig:2_4_31-Rel-Err-F7-C7b}
\par\end{centering}
\centering{}\vspace{0.6cm}
\includegraphics[scale=0.21]{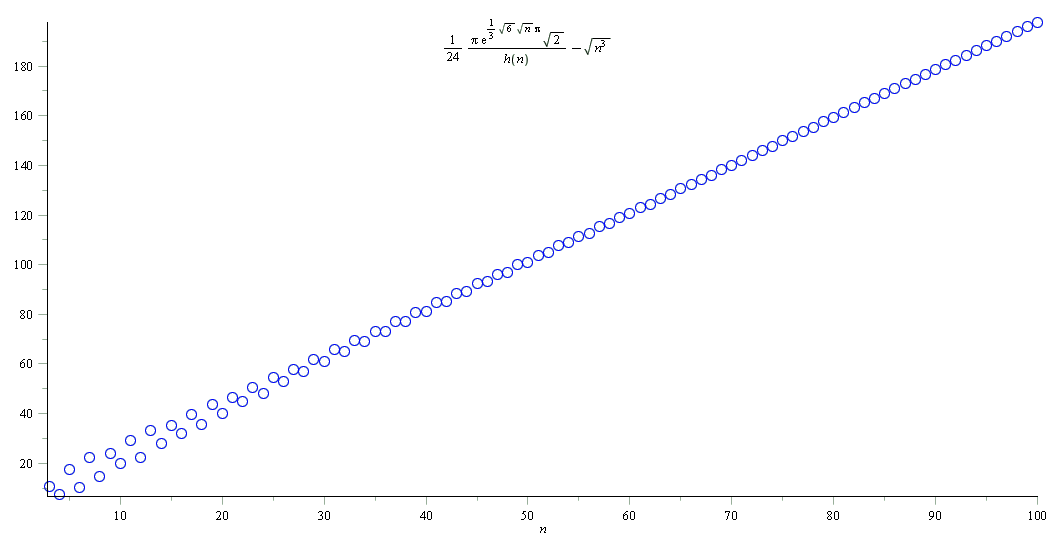}
\caption{The graph of the data $\left(n,\,C_{8}(n)\right)$}
 \label{Fig:2_4_32-(n,C8)}
\end{figure}

The relative error of 
\begin{equation}
F_{7a}(n)=I_{\mathrm{g}}(n)-\frac{\pi^{2}\exp\left(\sqrt{\frac{2}{3}}\pi\sqrt{n}\right)}{24\sqrt{3}n^{2}C_{7a}(n)}\label{eq:h(n)-estimation-3a}
\end{equation}
 and 
\begin{equation}
F_{7b}(n)=I_{\mathrm{g}}(n)-\frac{\pi^{2}\exp\left(\sqrt{\frac{2}{3}}\pi\sqrt{n}\right)}{24\sqrt{3}n^{2}C_{7b}(n)}\label{eq:h(n)-estimation-3b}
\end{equation}

to $h(n)$ when 1000 $\leqslant$ $n$ $\leqslant$ 10000 are shown
on Figure \ref{Fig:2_4_30-Rel-Err-F7-C7a} and Figure \ref{Fig:2_4_31-Rel-Err-F7-C7b}
(page \pageref{Fig:2_4_30-Rel-Err-F7-C7a}), respectively. In this
interval (1000, 10000), $F_{7a}(n)$ is obviously more accurate than
$F_{7b}(n).$ When $n$ $\leqslant$ 1000 the relative error of $\left\lfloor F_{7a}(n)+\frac{1}{2}\right\rfloor $
and $\left\lfloor F_{7b}(n)+\frac{1}{2}\right\rfloor $ are shown
on Table \ref{Table:Rel-Err-h(n)-F7a(n)} (page \pageref{Table:Rel-Err-h(n)-F7a(n)})
and Table \ref{Table:Rel-Err-h(n)-F7b(n)} ( page \pageref{Table:Rel-Err-h(n)-F7b(n)}).
In this case, $F_{7b}(n)$ is better than $F_{7a}(n)$. But neither
of them is as good as $I_{\mathrm{g1}}(n)$ or $I_{\mathrm{g2}}(n)$,
although they are more complicated than $I_{\mathrm{g1}}(n)$ and
$I_{\mathrm{g2}}(n)$.

\subsection{Estimate $h(n)$ When $n\leqslant100$\label{subsec:Estimate-less-100}}

All the estimation function for $h(n)$ found now are with very good
accuracy when $n$ is greater than 1000, but they are not so accurate
when $n<50$, especially when $n<25$. Although $I'_{\mathrm{g1}}(n)$
and $I'_{\mathrm{g2}}(n)$ are better than others, the relative error
are still greater than $1\permil$ when $n<40$.

\begin{table}[H]

\begin{centering}
\includegraphics[viewport=97bp 380bp 455bp 622bp,clip,scale=0.63]{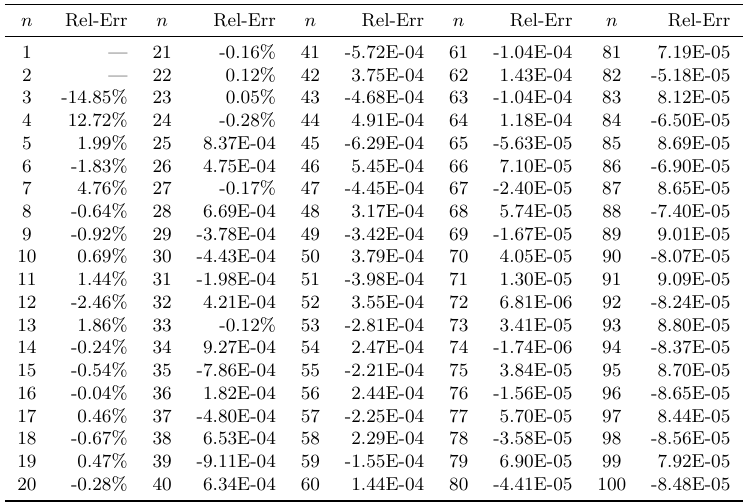}
\par\end{centering}
\caption{The relative error of $I_{\mathrm{g0}}(n)$ to $h(n)$ when $n\leqslant100$.}
 \label{Table:Rel-Err-h(n)-Ig0(n)}

\vspace{0.8cm}

\begin{centering}
\includegraphics[viewport=97bp 457bp 420bp 567bp,clip,scale=0.63]{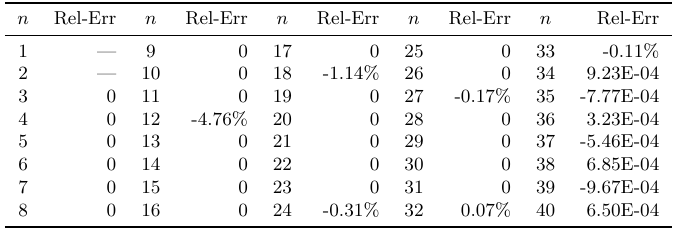}
\par\end{centering}
\caption{The relative error of $I'_{\mathrm{g0}}(n)$ to $h(n)$ when $n\leqslant100$.}
 \label{Table:Rel-Err-h(n)-Ig0(n)-round}

\end{table}

When $n<40$, it is too difficult to fit $\frac{\pi\exp\left(\pi\sqrt{\frac{2}{3}}\sqrt{n}\right)}{12\sqrt{2}h(n)}-n^{3/2}$
by a simple smooth function with high accuracy, as shown on Figure
\ref{Fig:2_4_32-(n,C8)} (on page \pageref{Fig:2_4_32-(n,C8)}). 
The figure of the points $\left(n,\ \frac{\pi\exp\left(\pi\sqrt{\frac{2}{3}}\sqrt{n}\right)}{12\sqrt{2}h(n)}-n^{3/2}\right)$
( $n$ = 3, 4, $\cdots$, 100) is not so complicated (as shown on
Figure \ref{Fig:2_4_32-(n,C8)}). It seems that we can fit them by
a simple piecewise function with 2 pieces, as the even points (where
$n$ is even) lie roughly on a smooth curve, so do the odd points.
If we try to fit them respectively, we will have the fitting function
below: 
\begin{equation}
C_{8}(x)=\begin{cases}
1.942141112\times x-0.4796781366\times\sqrt{x}\\
\quad+8.291226268,\quad n=3,5,7,\cdots,99; & \ \\
1.803056782\times x+2.356539877\times\sqrt{x}\\
\quad-6.043824511,\quad n=4,6,8\cdots,100. & \ 
\end{cases}\label{eq:C8(n)}
\end{equation}

Hence we can calculate $h(n)$ ( $3\leqslant n\leqslant100$) by 
\begin{equation}
h(n)\sim I_{\mathrm{g0}}(n)=\frac{\pi\exp\left(\sqrt{\frac{2}{3}}\pi\sqrt{n}\right)}{12\sqrt{2}\left(n^{3/2}+C_{8}(n)\right)},\quad3\leqslant n\leqslant100.\label{eq:h(n)-estimation-Ig0}
\end{equation}
 Consider that $h(n)$ is an integer, we can take the round approximation
of Equation \eqref{eq:h(n)-estimation-Ig0}, \label{Sym:I'g0(n)}
\nomenclature[I'g0(n)]{$I'_{\mathrm{g0}}(n)$}{The Ingham-Meinardus revised estimation formula when $ 3 \leqslant n \leqslant 100$. \pageref{Sym:I'g0(n)}}
\begin{equation}
I'_{\mathrm{g0}}(n)=\left\lfloor \frac{\pi\exp\left(\sqrt{\frac{2}{3}}\pi\sqrt{n}\right)}{12\sqrt{2}\left(n^{3/2}+C_{8}(n)\right)}+\frac{1}{2}\right\rfloor ,\quad3\leqslant n\leqslant100.\label{eq:h(n)-estimation-I'g0}
\end{equation}

Here $n$  begins from 3, not 1 or 2, because $\frac{I'_{\mathrm{g0}}(1)-h(1)}{h(1)}$
is meaningless since $h(1)=0$, and $I'_{\mathrm{g0}}(2)$ differs
from $h(2)$ a lot. Besides, the value of $h(1)$ and $h(2)$ are
clear by definition, so there is no need to use a complicated formula
to estimate them.

The relative error of $I_{\mathrm{g0}}(n)$ (or $I'_{\mathrm{g0}}(n)$)
to $h(n)$ are shown on Table \ref{Table:Rel-Err-h(n)-Ig0(n)} (or
Table \ref{Table:Rel-Err-h(n)-Ig0(n)-round}) on page \pageref{Table:Rel-Err-h(n)-Ig0(n)}.
Compared them with Table \ref{Table:Rel-Err-h(n)-Ig2(n)-round} on
page \pageref{Table:Rel-Err-h(n)-Ig2(n)-round}, we will find that
when $n\geqslant80$, $I'_{\mathrm{g2}}(n)$ is more accurate than
$I'_{\mathrm{g0}}(n)$; when $n<80$, $I'_{\mathrm{g0}}(n)$ is better.

\section{Conclusion}

 We have presented a recursion formula and several practical estimation
formulae with high accuracy to calculated the number $h(n)$ of conjugate
classes of derangements of order $n$, or the number of isotopy classes
of $2\times n$ Latin rectangles. 

If we want to obtain the accurate value of $h(n)$, we can use the
recursion formula \eqref{eq:hn-recursion} and write a program based
on it, while sometimes  we need to know the estimation value in a
program for technique reason especially when we use a general programming
language.

If we want to obtain the approximation value of $h(n)$ with high
accuracy, we can use the formulae \eqref{eq:h(n)-estimation-2}, \eqref{eq:h(n)-estimation-I'g0},
\eqref{eq:h(n)-estimation-1}, etc.

When $2\leqslant n\leqslant80$, we can use $I'_{\mathrm{g0}}(n)$
(Equation \eqref{eq:h(n)-estimation-I'g0}) , with a relative error
less than 0.11\% (while $32\leqslant n\leqslant80$) or mainly 0 with
very few exceptions (while $2\leqslant n\leqslant31$); when $n>80$,
we can use $I'_{\mathrm{g2}}(n)$ (Equation \eqref{eq:h(n)-estimation-2}).

When $n\geqslant100$, formulae $I'_{\mathrm{ga}}(n)$ (Equation \eqref{eq:Iga-round}),
$I'_{\mathrm{g1}}(n)$ (Equation \eqref{eq:h(n)-estimation-1}), $F_{7a}(n)$
(Equation \eqref{eq:h(n)-estimation-3a}) and $F_{7b}(n)$ (Equation
\eqref{eq:h(n)-estimation-3b}) are also very accurate although they
are not as good as Equations \eqref{eq:h(n)-estimation-2}.

With the asymptotic formula \eqref{eq:P_a,b} described in \cite{Daniel2006EleDerAsyPtFunc},
we can obtain some estimation formulae with high accuracy for some
other types of restricted partition numbers by the methods mentioned
in this paper or in \cite{liwenwei2016-Estmn-pn-arXiv}.

\section*{Data Availability}

The figures and tables used to support this study are included within
the article. The other data used to support this study are obtained
from a program made by the first author.

\section*{Conflicts of Interest}

The authors declare that there are no conflicts of interest regarding
the publication of this paper.

\section*{Acknowledgements}

This work was supported in part by the Fund of Research Team of Anhui
International Studies University, NO. awkytd1909.

The author would like to express the gratitude to the anonymous reviewers
for their valuable advice.

 \bibliographystyle{ieeetr}
\phantomsection\addcontentsline{toc}{section}{\refname}\bibliography{Ref-58}

\end{document}